\documentclass{siamart220329}
\usepackage{amsmath}
\usepackage{amsfonts}
\usepackage{amssymb}
\usepackage{dsfont}
\usepackage{bbm}
\usepackage{paralist}
\usepackage{xcolor}
\usepackage{enumitem}
\usepackage{optprog}

\usepackage{todonotes}
\newcommand\set[2]{\left\{ {#1}  :  {#2} \right\}}

\newcommand\R{\mathbb{R}}% reals
\newcommand\Z{\mathbb{Z}}% integers

\newcommand\x{\mathbf{x}}
\newcommand\y{\mathbf{y}}
\renewcommand\c{\mathbf{c}}
\newcommand\C{\mathbf{C}}
\renewcommand\u{\mathbf{u}}
\newcommand{\A}{\mathbf{A}}
\renewcommand{\b}{\mathbf{b}}
\newcommand{\E}{\mathbf{E}}
\newcommand{\B}{\mathbf{B}}
\newcommand{\Q}{\mathbf{Q}}
\newcommand{\f}{\mathbf{f}}

\newcommand\bz{\mathbf{0}}

\newcommand\s{\mathbf{s}}
\renewcommand\b{\mathbf{b}}

\renewcommand\int{\mathrm{int}\,}

\renewcommand\u{\mathbf{u}}

\newcommand{\1}{\mathbbm{1}}
\newcommand\abs[1]{\left| #1 \right|}

\newcommand\lcm{\textup{lcm}}
\newcommand\zIP{z_{IP}}

\newcommand\define{\,:=\,}
\newcommand{\Hu}{H_{\u}}
\newcommand{\Ho}{H_{\mathbf{0}}}
\newsiamremark{remark}{Remark}

\newcommand{\denom}{\textup{denom}}
\newcommand{\norm}[1]{\left|\left|#1\right|\right|}
\newcommand{\blambda}{\boldsymbol\lambda}
\newcommand{\bpi}{\boldsymbol\pi}

\newcommand{\calU}{U}
\newcommand\QPui{(\textbf{QP}^{\u}_{(i)})}
\newcommand\CPui{(\textbf{CP}^{\u}_{(i)})}
\newcommand\MICP{\textbf{(MICP)}}

\newcommand\ocalU{\overline{\calU}}

\newcommand\rec{\textup{\textbf{rec}}}
\headers{Exact ALD for MICP}{A. Bhardwaj, V. Narayanan, and A. Pathapati}
\title{Exact augmented lagrangian duality for mixed integer convex optimization}

\author{Avinash Bhardwaj\thanks{Industrial Engineering and Operations Research, Indian Institute of Technology Bombay, Mumbai, India 400076
  (\email{abhardwaj@iitb.ac.in}, \email{vishnu@iitb.ac.in}, \email{apathapati@iitb.ac.in}).}
\and Vishnu Narayanan\footnotemark[1]
\and Abhishek Pathapati\footnotemark[1]}

\begin{document}
\maketitle 
\begin{abstract}
Augmented Lagrangian dual augments the classical Lagrangian dual with a non-negative non-linear penalty function of the violation of the relaxed/dualized constraints in order to reduce the duality gap. We investigate the cases in which mixed integer convex optimization problems have an exact penalty representation using sharp augmenting functions (norms as augmenting penalty functions). We present a generalizable constructive proof technique for proving existence of exact penalty representations for mixed integer convex programs under specific conditions using the associated value functions. This generalizes the recent results for MILP (Feizollahi, Ahmed and Sun, 2017) and MIQP (Gu, Ahmed and Dey 2020) whilst also providing an alternative proof for the aforementioned along with quantification of the finite penalty parameter in these cases. 
\end{abstract}

% REQUIRED

% REQUIRED

\section{Introduction}
Given a polyhedral mixed-integer set $X \subseteq \Z^{n_1}\times\R^{n_2}$ and a real valued function $f:\R^n \mapsto \R$, consider the following mixed integer programming problem:
\begin{equation*}
    \textup{\textbf{(P)}} \qquad \zIP = \min\set{f(\x)}{\A\,\x\,= \b,\x\in X}
\end{equation*}
Solutions to mixed integer programming problems such as \textbf{(P)} are often computationally intractable, and the strong duality doesn't hold in general. As such, certain constraints of the optimization problem in context may be relaxed by using classical Lagrangian dual (LD) to yield good lower bounds on the optimal objective value. Specifically,
\begin{equation*}
    z^{LD} = \sup_{\blambda\in \R^m}\min_{\x\in X}(f(\x) + \blambda^{\top}(\b - \A\,\x)) \leq \zIP.
\end{equation*}
In contrast with the convex optimization problems, for nonconvex optimization problems such as \textbf{(P)}, classical Lagrangian dual may yield a non-zero duality gap, i.e. $z^{LD} < \zIP$. This duality gap may be avoided if the dual problem could be set up with, instead of affine dual functions, some other class of functions capable of penetrating possible ‘dents’ in the value function \cite{rockafellar2009variational}. Augmented Lagrangian dual (ALD), as the name suggests, augments the LD with a nonlinear penalty function of the violation of dualized constraints, 
\begin{equation*}
    z_{\rho}^{LD+} = \sup_{\blambda\in \R^m}\min_{\x\in X}(f(\x) + \blambda^{\top}(\b - \A\,\x) + \rho \psi(\b - \A\,\x))
\end{equation*}
where $\psi(.)>0$ is the augmenting penalty function and $\rho >0$ denotes the penalty parameter. Under certain conditions, a zero duality gap can be reached asymptotically by increasing the penalty parameter, $\rho$, to infinity \cite{wang2013nonlinear}. In some cases, the duality gap can be closed with a large enough finite value of the penalty parameter. In such cases, when the duality gap can be closed for a finite value of the penalty parameter $\rho$, the primal problem is termed to have an \textit{exact penalty representation}. Recently, the question of determining whether a class of non-convex optimization problems has an exact penalty representation has garnered quite some interest \cite{boland2015augmented,feizollahi2017exact,gu2020exact}. Boland and Eberhard \cite{boland2015augmented} use convex, monotone augmenting functions to close the mixed integer linear programming duality gap and also prove that for bounded pure-integer linear programs the gap can be closed for a finite penalty parameter. Additionally, the authors utilize techniques from \cite{burachik2005absence} to prove that duality gap is zero under the assumption that mixed integer linear programming problem attains the solution. Feizollahi et al. \cite{feizollahi2017exact} and Gu. et al \cite{gu2020exact} consider  level-bounded augmenting functions and prove that these augmenting functions close the duality gap as the penalty parameter goes to infinity for mixed integer linear programming (MILP) and mixed integer quadratic programming (MIQP) problems, respectively. They further prove that one can close the duality gap for a finite penalty parameter for sharp augmented Lagrangians. Burke \cite{burke1987exact,burke1991calmness} provides the characterization of the conditions when the duality gap can be closed for finite penalty parameter for augmented (sharp and proximal) Lagrangians. Another stream of research is focused on developing general approaches to solve augmented Lagrangian dual problems. Boland and Eberhard \cite{boland2015augmented} suggest the use of \textit{alternating directions method of multipliers} (ADMM) \cite{boyd2011distributed} for solving the MILP ALD problems. Cordova et al. \cite{cordova2021revisiting} provide a primal dual solution approach in form of a proximal bundle method to solve the ALD problems.\\

Value (Perturbation) functions of optimization problems provide a key insight into the properties of the augmenting functions that can be utilized for ALD to close the duality gap \cite{rockafellar2009variational,burke1987exact,burke1991calmness}. The structure of the value functions of MILPs has been extensively discussed in the literature \cite{meyer1975integer,blair1977value,blair1979value,blair1984constructive,blair1995closed,ralphs2014value}. Ralphs and Hassanzadeh \cite{ralphs2014value} provide an algorithm for the construction of the value function of a MILP and prove the finiteness of this algorithm in certain cases. It is also known that the value functions of both rational mixed-integer linear programs and continuous convex programs are lower semi-continuous \cite{meyer1975integer, bertsekas2003convex}.\\

The existence of exact penalty representations for MILPs \cite{boland2015augmented,feizollahi2017exact} and MIQPs \cite{gu2020exact} is well-established. However, the proof techniques used for MILPs and MIQPs are specific to these problem classes and aren't readily or necessarily generalizable. In the following, we present an alternative proof technique for proving existence of exact penalty representations in these cases using the associated value functions. This new proof technique helps us to generalize and prove the existence of exact penalty representations for mixed integer convex programs (MICPs) under specific conditions. Furthermore, this proof technique is also constructive in nature. Specifically, we provide an analytical form for the finite penalty parameter in case of MILPs and MIQPs and an upper bound in the case of MICPs.\\

The following discussion is organised in three sections. \Cref{sec:preliminaries} provides necessary definitions and highlights the notation used in the paper. \Cref{sec:results} outlines our primary results highlighting the developed generalizable contructive proof technique for proving existence of exact penalty representations in the case of MILPs, MIQPs and MICPs and construction of the penalty parameter in the respective cases. \Cref{sec:proofs} illustrates the proofs of the results presented in \Cref{sec:results}.

\section{Preliminaries}\label{sec:preliminaries}
Consider, to begin with, the following mixed integer convex programming (MICP) problem,
\begin{optprog*}
    {$\zIP = $ minimize}  & \objective{\hspace{3mm}f(\x)}\\
   \textbf{(MICP)} \qquad subject to & \A\,\x\,= \b\\
     & \x \in X \hspace{1mm}
\end{optprog*}
where $X:= \set{\x\in \Z^{n_1}\times\R^{n_2}}{\E\x\leq \f}$ is a polyhedral mixed integer set, $f:\R^n \mapsto \R$ is a real valued convex function, $\A,\E$ and $\b,\f$ are full rank matrices and vectors of appropriate dimensions, respectively. Through out the remainder of this discussion we assume that the matrices $\A, \E$ are rational matrices and $\b, \f$ are rational vectors and $f(\cdot)$ is differentiable. We further assume that both MICP and MICP's continuous relaxation are feasible and the corresponding optimal solutions exist.\\

Consider the Lagrangian relaxation of the MICP,
\begin{equation*}
    \textbf{(LR)} \hspace{1cm} z^{LR}(\blambda) = \min_{\x \in X} f(\x) + {\blambda}^{\top}(\b - \A\,\x).
\end{equation*}
The corresponding augmented Lagrangian relaxation and augmented Lagrangian dual (\cite{boland2015augmented,gu2020exact}) of the \MICP\, are defined as,
\begin{align*}
    \textbf{(ALR)} \hspace{0.5cm} z_{\rho}^{LR+}(\blambda) & = \min_{\x \in X} f(\x) + {\blambda}^{\top}(\b - \A\,\x) + \rho\, \psi(\b - \A\,\x)\\
    \textbf{(ALD)} \hspace{1cm} z_{\rho}^{LD+} & = \sup_{\blambda\in \R^m}\min_{\x \in X} f(\x) + {\blambda}^{\top}(\b - \A\,\x) + \rho\, \psi(\b - \A\,\x)
\end{align*}
where $\psi:\R^m \mapsto \R$ is a real valued function. We further designate $\psi$ to have an \textit{exact penalty representation} if $\exists\, 0 < \rho < \infty$, such that $z^{LD+}_{\rho} = \zIP$. In the case $\psi(\cdot) = \norm{\cdot}$ we augmented Lagrangian relaxation is referred to as sharp augmented Lagrangian relaxation. In particular, 
$$
    \textbf{(SALR)} \hspace{1cm} z_{SALR}(\rho) = \min_{\x \in X} f(\x) + \rho\, \norm{\b - \A\,\x}\\
$$
It follows that $z^{LR+}_{\rho}(\blambda) \leq \zIP$ and $z^{LR+}_{\rho}(\bz) = z_{SALR}(\rho)$. 
\subsection{Notation}\label{sec:notation} The sets $\R^n,\Z^n,\mathbb{Q}^n$ denote the set of real numbers, integers and rational numbers in $n-$dimensional vector space, respectively. Additionally, $\R^n_+,\Z^n_+,\mathbb{Q}^n_+$ denote the non-negative counterparts of the respective sets. We use $[n]$ to denote the index set $\{1,2 \ldots, n\}$. For an index set $E \subset \Z^+$ and $y\in \Z$, we define the translation of an index set as $y + E \define \set{x + y}{x \in E}$. For a mixed integer set $S \subseteq \R^n \times \Z^m$, we denote by $S_R$ the continuous relaxation of $S$. Throughout the following discussion we define $I$ and $C$ as the index sets corresponding to integer and continuous variables, respectively. Analogously, for any $\u \in \R^n$, let $\u = (\u_I,\u_C)$ denote the partition of $\u$ into the integer and continuous variables, respectively. Similarly for a symmetric matrix $\Q \in \R^{n\times n}$, let $Q_{UV}$ denote the submatrix of $\Q$ formed by rows indexed by $U$ and columns indexed by $V$. For $\delta > 0$ and $\bar{\x} \in \R^n$, we define by $\mathcal{N}_{\delta}(\bar{\x}) \define \set{\x \in \R^n}{\norm{\x -\bar{\x}}\leq \delta}$, the ball with center $\bar{\x}$ and radius $\delta$. The recession cone of a set $S$ is denoted by $\rec(S)$. We also define, $\rec(f)$, \emph{recession cone of a convex function} $f:\R^n \mapsto \R$ as the recession cone of a level set $\set{x \in \R^n}{f(\x)\leq \alpha}$, for some $\alpha \in \R$ (\cite{bertsekas2003convex}, Proposition 2.3.1). For given convex sets $A$ and $B$, we say that the sets $A$ and $B$ have \emph{no common/distinct non-zero directions of recession} if $\rec(A)\cap \rec(B)\setminus \{\bz\} = \emptyset$.\\

We say that a function $g : \R^m \mapsto \R$ is $L$ smooth, or alternatively, has \textbf{Lipschitz-continuous} gradients if there exists a Lipschitz
constant, $L < \infty$, such that
$$
    \norm{\nabla{g(\x)} - \nabla{g(y)}}\leq
    L\norm{\x-\y}~~~~ \forall~~ \x,\y \in \R^m.
$$
Additionally, we say a function $g : \R^m \mapsto \R$ is $\mu$ \textbf{strongly-convex} if there exists $ 0 <\mu < \infty$, such that
$$
    \norm{\nabla{g(\x)} - \nabla{g(y)}}\geq\mu\norm{\x-\y}~~~~ \forall~~ \x,\y \in \R^m.
$$
$\psi : \R^m \mapsto \R_{+}$ is termed as a proper, non-negative, lower semi-continuous function, level-bounded augmenting function \emph{iff} 
$$
\psi(\bz) = 0,~ \psi(\u)> 0 ~~\forall~~ \u \neq \bz, ~~\text{and}~~\textbf{diam}\set{\u}{\psi(\u)\leq \delta} < +\infty ~~\forall~~ \delta > 0.
$$ 
Moreover, $\lim_{\delta\downarrow 0}\textbf{diam}\set{\u}{\psi(\u)\leq \delta} = 0$ \cite{feizollahi2017exact}.\\

Consider the hyperplane $H_{\u} \define \set{\x\in\R^n}{\A\x = \b + \u}$, and let the corresponding value function associated with {\MICP} be defined as 
\begin{align}
    \phi(\u) = \min\set{f(\x)}{\x \in X\cap\Hu}.\label{def:valuefunction}
\end{align}
In particular, $\zIP = \phi(\bz) = \min_{\x\in X \cap \Ho} f(\x)$. Finally, we define the set $\calU$ as the set of all possible perturbation vectors $\u$ such that the feasible set of \cref{def:valuefunction} is non-empty. In particular, $\calU \define \set{\u\in \R^m}{X \cap \Hu \neq \emptyset}$.

\section{Main results}\label{sec:results}
The primary contribution of this work entails a generalizable constructive proof technique for proving existence of exact penalty representations for mixed integer convex programs under specific conditions. We would like to emphasize, in particular, that while the results on existence of exact penalty represenatations using sharp Lagrangians in the specific case of rational MILPs and bounded integer variable MIQPs have been discussed in literature \cite{feizollahi2017exact,gu2020exact}, the proofs use specific properties of MILPs and MIQPs and thus don't necessarily generalize to MICPs. The proposed proof technique utilizes the properties of the associated value functions. Specifically, the proof utilizes \cref{lem:3} in conjunction with the lower semi-continuity of the value functions of continuous convex optimization problems in both the aforementioned cases. In addition to resolving the existence of an exact penalty representation in the cases discussed, we further provide a quantification of the associated penalty parameter $\rho$. This quantification, to the best of our knowledge, has not been discussed in literature.\\

The following theorems formalize this discussion. The proofs of the theorems follow in \cref{sec:proofs}. It should be noted that in addition to the assumptions stated in \cref{sec:preliminaries}, we further assume that $\c = (\c_C,\c_I)$ is a rational vector and $\Q$ is a rational symmetric positive semi-definite matrix.
\begin{theorem}\label{thm:1}
Consider the following mixed integer linear programming problem,
\begin{optprog*}
    minimize  & \objective{\hspace{5mm}\c_I^{\top}\x_I + \c_C^{\top}\x_C}\\
    \textup{\textbf{(MILP)}}\hspace{2mm} subject to & \A_I\,\x_I + \A_C\x_C = \b\hspace{6.5mm}\\
     & \E_I\x_I + \E_C\x_C \leq \f\hspace{7mm}\\
     & (\x_I, \x_C) \in \Z^{n_1}\times\R^{n_2}\hspace{0.3cm}
\end{optprog*}
There exists an exact penalty representation for \textup{\textbf{(MILP)}}. Furthermore, the finite penalty parameter $\rho$ depends on $\A_C$ and $\c_C$ and $\c_I$ .
\end{theorem}

\begin{theorem}\label{thm:2}
Consider the following mixed integer quadratic programming problem,
\begin{optprog*}
    minimize  & \objective{\frac{1}{2}\x_C^{\top}Q_{CC}\x_C + \frac{1}{2}\x_I^{\top}Q_{II}\x_I + \x_I^{\top}Q_{IC}\x_C- {\c}_I^\top\x_I - {\c}_C^\top\x_C}\\
    \textup{\textbf{(MIQP)}} \hspace{2mm}  subject to & \hspace{1.2cm} \A_I\,\x_I + \A_C\x_C = \b\hspace{6.5mm}\\
    & \E_I\x_I + \E_C\x_C \leq \f \hspace{7mm}\\
    & \norm{\x_I}_{\infty}\leq M \hspace{0.5cm}\\
    & (\x_I, \x_C) \in \Z^{n_1}\times\R^{n_2}\hspace{0.3cm}
\end{optprog*}
There exists an exact penalty representation for \textup{\textbf{(MIQP)}}. Furthermore, the finite penalty parameter $\rho$ depends on  $\A,\Q,{\c}$ and $M$.
\end{theorem}
Generalizing the discussion on MILPs and MIQPs our next result primarily focusses the discussion on exact penalty representations of MICP's under specific conditions. We initiate the discussion by proving that the Pure Integer Convex Programs (PICPs) with rational data have an exact penalty representation. We further prove that MICPs where either the objective function is strongly convex or where the recession cone of the epigraph of the objective function and the recession cone of the continuous relaxation of the feasible set have no common non-zero directions of recession the duality gap can be closed with level-bounded functions asymptotically as $\rho \to \infty$. As our concluding result, we establish in \cref{thm:3} that (a) MICPs where the recession cone of the epigraph of the objective function and the recession cone of the continuous relaxation of the feasible set have no common non-zero directions of recession and (b) MICPs with bounded integer variables, have exact penalty representation when using norms as augmenting functions.

\begin{theorem}\label{thm:3}
Consider the following mixed integer convex programming problem,
\begin{optprog*}
    minimize  & \objective{\hspace{1cm}f({\x_I},{\x_C})}\\
    \textup{\MICP} \hspace{2mm} subject to & \hspace{0cm} \A_I\,\x_I + \A_C\x_C = \b\hspace{6.5mm}\\
    & \E_I\x_I + \E_C\x_C \leq \f \hspace{7mm}\\
    % & \norm{\x_I}_{\infty}\leq M \hspace{0.5cm}\\
    % & \x_I\in \Z_{+}^{n_1}\hspace{0.4cm}\\ 
     & (\x_I, \x_C) \in \Z^{n_1}\times\R^{n_2}\hspace{0.3cm}
\end{optprog*}
Let $F = X\cap \Ho$ denote the feasible region of \textup{\MICP}.\\
\begin{enumerate}[label = (\alph*)]
    \item If $\norm{\x_I}_{\infty}\leq M,$ $\forall\,\x \in F$ then there exists an exact penalty representation for \textup{\MICP}.
    \item If the recession cone of $f$ and recession cone of $F_R$, the feasible set of continuous relaxation of \textup{\MICP}, have no common non-zero directions of recession, then there exists an exact penalty representation for \textup{\MICP}.
\end{enumerate} 
Additionally, if $f$ is $\mu$-strongly convex and $L$ smooth, then the finite penalty parameter $\rho = \mathcal{O}\left(\frac{L\beta\gamma}{\mu}\right)$ where $\displaystyle \beta \define \max_{B \in \mathcal{B}}\norm{B^{-1}}_F$, $\mathcal{B}$ being the set of all possible invertible submatrices of $[\A^T_C \,\, -\A^T_C \hspace{2mm}-{\E}^T_C]$ and $\gamma$ can be explicitly computed given $f(\bz)$, $\norm{\nabla{f(\bz)}}$ and $f(\bar{\x})$ for any $\bar{\x} \in F$.
% $\gamma = \sqrt{f(\bar{\x})}$ for any $\bar{\x} \in F$.\\
\end{theorem}
\begin{remark}
The parameter $\gamma$ of \cref{thm:3}(c) can be explicitly computed as $\gamma = 2\norm{\nabla{f(\bz)}} + \sqrt{\norm{\nabla{f(\bz)}}^2 + 2\mu\,\left(f(\bar{\x})- f(\bz)\right)}$. Observe that the feasible region $F$ of \MICP{} is a rational polyhedron. From Corollary 17.1d of~\cite{schrijver1986theory}, there exists a feasible $\bar{\x} \in F$ whose size is polynomially bounded by the size of $\A, \b, \E$, and $\f$. Therefore, the parameter $\gamma$ of \cref{thm:3} is ``small'' with respect to the input size if the function value $f(\cdot)$, and the gradients of $f(\cdot)$ can be computed efficiently. 
\end{remark}
\section{Proofs of the Theorems}\label{sec:proofs}
\subsection{Preliminary Lemmas}
\begin{lemma}\label{lem:1}
$\displaystyle z_{SALR}(\rho) = \min_{\u \in \calU}\,\, \phi(\u) + \rho\norm{\u}$ . 
\end{lemma}
\begin{proof}
Observe the following inequality which follows immediately from the aforementioned definitions
\begin{align}
    z_{SALR}(\rho) = \min_{\x \in X} f(\x) + \rho\, \norm{\b - \A\,\x} & \leq \min_{\x\in X \cap \Hu} f(\x) + \rho\norm{\u} ,\,\,\,\,\, \forall\,\u \in \calU \nonumber\\
                    & = \phi(\u) + \rho\norm{\u} ,\,\,\,\,\, \forall\,\u \in \calU. \label{eq:1}
\end{align}
Conversely, it holds that $\forall\,\u \in \calU$
\begin{align*}
    \phi(\u) + \rho\norm{\u} & = \min_{\x\in X \cap \Hu} f(\x) + \rho\norm{\u}\\
                             & \leq f(\x) + \rho\norm{\A\x-\b} ,\,\,\,\,\, \forall\,\x\in X\cap \Hu.
\end{align*}
This implies that $ \min_{\u \in \calU} \phi(\u) + \rho \norm{\u} \leq f(\x) + \rho \norm{\A\x - \b} $ $\forall~ x\in X \cap \Hu,\,\forall~\u \in U$ which yields that $ \min_{\u \in \calU} \phi(\u) + \rho \norm{\u} \leq f(\x) + \rho \norm{\A\x - \b} $ $\forall~ \x \in X $. This implies,
\begin{align}
    \min_{\u\in \calU} \phi(\u) + \rho\norm{\u} \leq \min_{\x\in X} f(\x) + \rho\norm{\A\x-\b}. \label{eq:2}
\end{align}
The result follows from \cref{eq:1} and \cref{eq:2}.
\end{proof}
Consider the partition of the feasible set $X = F_> \cup F_{\leq}, F_> \cap F_{\leq} = \emptyset$ where $F_{\leq} = \set{\x \in X}{f(\x) \leq \zIP}$ and $F_> =\set{\x \in X}{f(\x) > \zIP}$. Observe that for $\x \in F_>$ any positive $\rho > 0$ implies $f(\x) + \rho\norm{\b - \A\,\x} > \zIP$. 
Thus, for MICP to have an exact penalty representation, it suffices to show that $\exists$ $0 < \rho < \infty$ such that $f(\x) + \rho \norm{\b - \A\,\x} \geq \zIP$, $\forall \x \in F_{\leq}$.

\begin{proposition}\label{prop:2}
Consider the continuous relaxation of MICP, and let the optimal objective value $z_R = \min_{X_R\cap \Ho} f(\x)$ be attained at $\x_R$. 
Let $\blambda_{\A}$ and $\blambda_{\E}$ be the Lagrange multipliers for the constraints $\A\x = \b$ and $\E\x\leq \f$, respectively.\vspace{1mm}
\begin{enumerate}[label = \roman*)]
    \item $\phi(\u) > -\infty$ for all $\u \in \calU$.\vspace{1mm}
    \item If $\exists ~ \boldsymbol{\alpha} \in X_R$ such that $f(\boldsymbol{\alpha}) + \rho\norm{\b-\A\,\boldsymbol{\alpha}}\leq \zIP$ and $\rho > \norm{\blambda_\A}$  then 
        \[
                \norm{\A\,\boldsymbol{\alpha} - \b} \leq \frac{\zIP-z_R}{\rho -\norm{\blambda_{\A}}}.
        \]    
\end{enumerate}
\end{proposition}
\begin{proof}
Consider the Lagrangian function for the continuous relaxation of MICP, i.e. for $\x \in \R^n$
\begin{equation*}
    \mathcal{L}(\x,\blambda_{\A},\blambda_{\E}) = f(\x) - \blambda_\A^{\top}(\A\,\x-\b)-\blambda_\E^{\top}(\f-\E\x)
\end{equation*}
As the relaxed problem has only affine constraints strong duality holds (Chapter 5,Section 5.2.3\cite{boyd2004convex}).

 As strong duality holds, the first order necessary (KKT) conditions for $\mathcal{L}(\x,\blambda_{\A},\blambda_{\E})$ can be characterized as
\begin{align}\label{eq:4}
\begin{split}
    \nabla{f(\x)} = \A^{\top}\blambda_\A - \E^{\top}{\blambda_\E}\\
    \blambda_\E^{\top}(\f - \E\x) = 0\\
    \blambda_\A^{\top}(\A\,\x - \b) = 0\\
    \blambda_\E \geq \bz
\end{split}
\end{align}
Observe that $\x_R$ satisfies the system of equations \cref{eq:4}.\\
\begin{enumerate}[label = \roman*)]
    \item To see the result, observe that it suffices to show that $\phi(\u) > -\infty$ for all $\u \in \calU$ such that $\phi(\u) \leq \phi(\bz) = \zIP$. Further, consider the set $F_{\leq}$ defined as,
            $$
                F_{\leq} \define \set{\x \in X}{f(\x) \leq \zIP}.
            $$
            We have for all $\x \in F_{\leq}$, $\zIP \geq f(\x) \geq f(\x_R) + \nabla{f(\x_R)}(\x-\x_R)= z_R + \nabla{f(\x_R)}(\x-\x_R)$. Substituting from \cref{ineq:gradient} we obtain, for all $\x \in F_{\leq}$,
            $$
                \zIP \geq f(\x) \geq z_R + \blambda_{\A}^{\top} (\A\,\x - \b) -  \blambda_\E^{\top} (\E\x -  \f )
            $$
            As $\blambda_{\E} \geq \bz$ and $\x \in X_R$,
            \begin{align*}
                \zIP \geq f(\x) & \geq z_R + \blambda_{\A}^{\top} (\A\,\x - \b)\\
                & \geq z_R - \norm{\blambda_{\A}}\norm{(\A\,\x - \b)}\\
                & = z_R - \norm{\blambda_{\A}}\norm{\u} > -\infty.
            \end{align*}
            Since $f(\x) > -\infty$ for all $\x \in F_{\leq}$, the result follows.\\
    \item Since we have for $\boldsymbol{\alpha}\in X_R$, $f(\boldsymbol{\alpha}) + \rho\norm{\A\,\boldsymbol{\alpha} - \b} \leq  \zIP$; convexity of $f$ yields,
            \begin{align}
                f(\x_R)+ \nabla{f(\x_R)}^{\top}(\boldsymbol{\alpha}-\x_R) +\rho \norm{\A\,\boldsymbol{\alpha} - \b} &\leq \zIP , \text{ i.e.,}\nonumber \\
                z_R+ \nabla{f(\x_R)}^{\top}(\boldsymbol{\alpha}-\x_R) +\rho \norm{\A\,\boldsymbol{\alpha} - \b} &\leq \zIP. \label{ineq:1}
            \end{align}
            It follows from strong duality that
            \begin{align}
            \nabla{f(\x_R)}^{\top}(\x-\x_R) &= \blambda_{\A}^{\top} \A(\x - \x_R) -\blambda_\E^{\top} \E(\x-\x_R)\nonumber\\
             &= \blambda_{\A}^{\top} (\A\,\x - \b) -  \blambda_\E^{\top} (\E\x -  \f ). \label{ineq:gradient}
            \end{align}
            Substituting in \cref{ineq:1}, we obtain
            \begin{align*}
            z_R+ \blambda_{\A}^{\top} (\A\,\boldsymbol{\alpha} - \b) +  \blambda_\E^{\top} (\f - \E\boldsymbol{\alpha}) +\rho \norm{\A\,\boldsymbol{\alpha} - \b}  & \leq \zIP .
            \end{align*}
            As $\blambda_{\E} \geq \bz$ and $\boldsymbol{\alpha} \in X_R$, the above inequality can be rewritten as
            % \blambda_{\A}^{\top} (\A\,\x - \b) +   \blambda_\E^{\top} ( \f - \E\x ) +\rho \norm{\A\,\x - \b} & \leq \zIP - z_R\\
            \begin{align}
            \blambda_{\A}^{\top}(\A\,\boldsymbol{\alpha} - \b) + \rho \norm{\A\,\boldsymbol{\alpha} - \b} & \leq \zIP - z_R \label{ineq:2}
            \end{align}
            Using Cauchy Schwarz inequality in \cref{ineq:2} yields for all $\x \in X_R$,
            \begin{align*}
                -\norm{\blambda_\A}\,\norm{\A\,\boldsymbol{\alpha} - \b} + \rho\,\norm{\A\,\boldsymbol{\alpha} - \b} \leq \zIP - z_R\\
            (\rho - \norm{\blambda_\A})\norm{\A\,\boldsymbol{\alpha} - \b}\leq \zIP - z_R\\
            \norm{\A\,\boldsymbol{\alpha} - \b} \leq \frac{\zIP-z_R}{\rho - \norm{\blambda_\A}}.
            \end{align*}
\end{enumerate}
\end{proof}
%%%%%%%%%%%%%%%%%%%%%%%%%%%%%%%%%%%%%%%%%%%%%%%%%%%%%%%%%%%%%%%%%%%%%%%%%%%%%%%%%%%%%%%%%%%%%%%%%%%%%%%%%%%%%%%%%%%%%%%%%%%%%%%%%%%%%%%%%%%%%%
\begin{corollary}\label{cor:3}
Consider the set $U_\rho = \set{\u \in \calU}{\phi(\u) + \rho\norm{\u} \leq \phi(\bz)}$. If $\rho > \norm{\blambda_{\A}}$ then for all $\u \in U_\rho$, $$\norm{\u} \leq \frac{\zIP-z_R}{\rho -\norm{\blambda_\A}}.$$
\end{corollary}
\begin{proof}
If for any $~ \u \in U_{\rho},~\exists~ \x \in X\cap H_{\u}$ such that $f(\x) + \rho\norm{\A\x - \b} \leq z_{IP}$, the result follows from \Cref{prop:2}. Alternatively, if $f(\x) + \rho\norm{\A\x - \b} > z_{IP}\,\, \forall~ \x \in X\cap H_{\u}$ for some $\u \in U_{\rho}$. 
\begin{align*}
\phi(\u) + \rho(\norm{\u}) =  \inf_{X\cap H_{\u}} f(\x) + \rho\norm{\A\x - \b} \geq z_{IP}
\end{align*} 
If \(\inf_{X\cap H_{\u}} f(\x) + \rho\norm{\A\x - \b} > z_{IP}\) then $\u \not\in U_{\rho}$ leading to a contradiction. If \(\inf_{X\cap H_{\u}} f(\x) + \rho\norm{\A\x - \b} = z_{IP}\) and there exists an \(\x \in X\cap H_{\u}\) such that \(f(\x) + \rho \norm{\A\x - \b} = z_{IP}\) then the 
result follows from \cref{prop:2}. If there doesn't exist \(\x \in X\cap H_{\u}\) such that \(f(\x) + \rho\norm{\A\x - \b}  = z_{IP}\) then one can find a sequence of \(\x_{p}\) such that the follows holds
\begin{align*}
    f(\x_{p}) + \rho \norm{\A\x_p - \b} \leq z_{IP} + \frac{1}{p} \text{ (From definition of infimum)}\\
\end{align*}
Let \(\x_R\) be the optimal solution of the optimization program \(\inf_{\x \in F_R} f(\x)\). Since we are minimizing a convex function over a rational polyhedron the KKT conditions hold as follows.
\begin{align*}   
\nabla{f(\x_{R})} = \A^{\top}\blambda_\A - \E^{\top}{\blambda_\E}\\
    \blambda_\E^{\top}(\f - \E\x_R) = 0\\
    \blambda_\A^{\top}(\A\,\x_R - \b) = 0\\
    \blambda_\E \geq \bz
\end{align*}
\begin{align*} 
    f(\x_R)+ \nabla{f(\x_R)}^{\top}(\x_{p}-\x_R) +\rho \norm{\A\,\x_p - \b} &\leq \zIP + \frac{1}{p}\text{ (Convexity)}\\
    z_{R} + \blambda_{\A}^{\top} (\A\,\x_p - \b) +  \blambda_\E^{\top} (\f - \E\x_p) +\rho \norm{\A\,\x_p - \b} &\leq \zIP + \frac{1}{p}\\
    % (\rho - \norm{\blambda_{\A}}) \norm{\u} &\leq \frac{\zIP - z_{R}}{\rho - \norm{\blambda_\A}} + \frac{1}{p(\rho - \norm{\blambda_{\A}})}\\
    \norm{\u} &\leq \frac{\zIP - z_{R}}{\rho - \norm{\blambda_{\A}}} + \frac{1}{p(\rho - \norm{\blambda_{\A}})}
\end{align*}
where the last inequality follows from non-negativity of $\blambda_\E$ and $(\f - \E\x_p).$ As \(p\) is arbitrary, we have as $p \to \infty$,
\begin{align*}
    \norm{\u} &\leq \frac{\zIP - z_{R}}{\rho - \norm{\blambda_{\A}}}
\end{align*}
\end{proof} 
\begin{corollary}\label{cor:4}
Consider the set $\,\ocalU \define \set{\u \in \calU}{\phi(\bz) \geq \phi(\u)}$. For $\alpha \geq 0$, $\phi(\bz) \leq \phi(\u) + \alpha$ $\forall \,\u \in \calU$ if $\phi(\bz) \leq \phi(\u) + \alpha$ $\forall \,\u \in \ocalU$.
\end{corollary}
\begin{proof}
    For $\u \in \calU\setminus\ocalU$, $\phi(\bz) < \phi(\u) \leq \phi(\u) + \alpha$. The result follows.
\end{proof}
\begin{lemma}\label{lem:3}
If there exist $\delta, \kappa > 0$ such that  $\phi(\bz) \leq \phi(\u) + \kappa \norm{\u}$ for all $\u \in \mathcal{N}_{\delta}(\bz) \cap \calU$,  then there exists $0 < \rho^{*} < \infty$ such that $z_{SALR}(\rho) = \zIP$ for all $\rho > \rho^{*}$.
\end{lemma}
\begin{proof}
Observe that, for $\rho > 0$ it suffices to consider $\u \in U_\rho$. We have, from \cref{lem:1} 
$$ z_{SALR}(\rho) = \min_{\u\in \calU} \phi(\u) + \rho \norm{\u}.$$
It follows that,
\begin{align*}
    z_{SALR}(\rho) = \min_{\u\in U_\rho}\phi(\u) + \rho \norm{\u}
\end{align*}
\cref{prop:2} yields that for $\rho > \norm{\blambda_{\A}}$ and $\u \in U_{\rho}$, $\norm{\u}$ is bounded. Thus, we have
\begin{equation*}
    z_{SALR}(\rho) = \min_{\u \in \mathcal{N}_{\delta_{\rho}}(\bz) \cap\, U_\rho}\phi(\u) + \rho \norm{\u} 
\end{equation*}
where $\displaystyle\delta_{\rho} = \frac{\zIP - z_R}{\rho -\norm{\blambda_{\A}}}$. From \cref{prop:2} we see that as $\rho \to \infty,\, \mathcal{N}_{\delta_{\rho}}(\bz) \cap U_\rho \to \{\bz\}$. If $\exists$ $\delta > 0$, and $0 < \kappa < \infty$ such that $\phi(\bz) \leq \phi(\u) + \kappa \norm{\u}$ for all $\u \in \mathcal{N}_{\delta}(\bz) \cap U_\rho$, it follows that,
\begin{align}
    \min_{\u\in \mathcal{N}_{\delta}(\bz)\cap\, U_\rho} \phi(\u) + \kappa \norm{\u} = \phi(\bz) = \zIP\label{eq:9}
\end{align}
Alternatively, we can increase $\rho$ such that $\forall\, \u \in U_\rho$, $\norm{\u} \leq \dfrac{\zIP - z_R}{\rho -\norm{\blambda_{\A}}} < \delta$.  Rearranging the terms, we obtain $\rho > \norm{\blambda_\A} + \frac{\zIP - z_R}{\delta}$.  It follows from \cref{eq:9} that for any choice of $\displaystyle \rho > \max\left\{\kappa, \norm{\blambda_\A} + \frac{\zIP - z_R}{\delta}\right\}$, we have  
\begin{align*}
    z_{SALR}(\rho) = \min_{\u\in \mathcal{N}_{\delta_{\rho}}(\bz) \cap\, U_\rho} \phi(\u) + \rho \norm{\u} = \zIP.
\end{align*}
\end{proof}
\subsection{Proof of \cref{thm:1}}
\begin{proposition}\label{prop:4}
Consider the value function $\phi(\u)$ of a mixed integer linear programming problem,
\begin{optprog*}
    {$\phi(\u) = $ minimize}  & \objective{\hspace{5mm}\c_I^{\top}\x_I + \c_C^{\top}\x_C}\\
    \hspace{2mm} subject to & \A_I\,\x_I + \A_C\x_C = \b + \u\\
     & \x_I\in \Z_{+}^{n_1}\hspace{1.2cm}\\ 
     & \x_C \in \R_{+}^{n_2}\hspace{1.2cm}
\end{optprog*}
There exists a $\delta > 0$ such that for every $ \u\in \mathcal{N}_{\delta}(\bz)\cap \calU$, $ \phi(\bz) \leq \phi(\u) + \Gamma\norm{\u}$ 
where $\Gamma$ is a constant which depends on $\A_C$ and $\c_C$ and $\c_I$ .
\end{proposition} 
\begin{proof}
From \cref{cor:4}, it suffices to show the existence of a $\delta > 0$ such that for every $ \u\in \mathcal{N}_{\delta}(\bz)\cap \ocalU$, $ \phi(\bz) \leq \phi(\u) + \Gamma\norm{\u}$. Observe that, if there doesn't exist a limiting sequence to $\bz$ in $\ocalU$ then there exists a $\delta > 0$ such that $\ocalU\cap \mathcal{N}_{\delta}(\bz) = \{\bz\}$. Hence $\phi(\bz)\leq \phi(\u)$ for all $\u \in \mathcal{N}_{\delta}(\bz)$. The result follows in this case.

Conversely, if there does exist a limiting sequence to $\bz$ in $\ocalU$ then given $\phi(\u) \leq \phi(\bz),\,\forall\,\u\in \ocalU$ implies that $\displaystyle \limsup_{\u \rightarrow \bz}\phi(\u) \leq \phi(\bz)$.\\ 
Since $\phi(\u)$ is lower semi-continuous\cite{meyer1975integer}, we have $\displaystyle\liminf_{\u\rightarrow \bz} \phi(\u) \geq  \phi(\bz)$. It follows that $\displaystyle\lim_{\u \rightarrow \bz} \phi(\u) = \phi(\bz)$.

For a given $\u \in \ocalU$, let $\displaystyle \x^\u = (\x^\u_I,\x^\u_C) = \arg\min_{X\cap \Hu} \phi(\u)$. Rewriting $\x^\u_C$ in terms of $\x^\u_I$ we have 
\begin{equation}\label{eq:11}
    \x^\u_C = \A_{\B_{\u}}^{-1}(\b + \u - \A_I\x^\u_I) \text{  and  } \phi(\u) = \c^{\top}_I\x^\u_I + \c_{\B_{\u}}^\top\A_{\B_{\u}}^{-1}(\b + \u - \A_I\x^\u_I),
\end{equation}
where $\B_{\u}$ is the optimal basis with respect to the continuous part of the solution, and $\A_{\B_{\u}}$ and $\c_{\B_\u}$ are the columns of $\A_C$ and elements of $\c_C$ indexed by $\B_{\u}$. Now, for every $\epsilon > 0$ there exists a $\delta > 0$ such that
\begin{equation*}
    % \abs{\phi(\u) - \phi(\bz)} < \epsilon \hspace{1cm} \forall \norm{\u} < \delta.
        % \abs{\phi(\u) - \phi(\bz)} < \epsilon \hspace{1cm} \forall \mathcal{N}_{\delta}(\mathbf{0})\cap \overline{U}
       \abs{\phi(\u) - \phi(\bz)} < \epsilon \hspace{1cm} \forall~ \u \in \mathcal{N}_{\delta}(\mathbf{0})\cap \overline{U}
\end{equation*}
Substituting for $\phi(\u)$ from \cref{eq:11},
\begin{align*}
    \abs{\c_I^\top\x^\u_I + \c_{\B_\u}^\top\A_{\B_\u}^{-1}(\b+ \u - \A_I\x^\u_I) - \phi(\bz)} < \epsilon 
\end{align*}
\begin{equation}\label{eq:12}
    \abs{\c_I^\top\x^\u_I + \c_{\B_\u}^\top\A_{\B_\u}^{-1}(\b - \A_I\x^\u_I) -\phi(\bz)+ \c_{\B_\u}^\top\A_{\B_\u}^{-1}\u} < \epsilon
\end{equation}
Since $\A$ and $\b$ are rational, one can assume without loss of generality that $\A$ and $\b$ are integral. This implies that $\A_{\B_\u}$ is integral. Consequently, $\det(\A_{\B_\u})$ and $\text{Adj}(\A_{\B_\u})$ are integral as well, where $\det(\mathbf{M})$ and $\text{Adj}(\mathbf{M})$ denote the determinant and adjugate matrix of $\mathbf{M}$, respectively. Consider the function $\denom: \mathbb{Q}\to \Z$,
\begin{align*}
    \denom(r) =  \begin{cases} 
      \abs{q} & r = \dfrac{p}{q}\,\, \text{such that   } \textup{gcd}(p,q) = 1\,\, r\neq 0, \\[2mm]
      1 & r = 0
   \end{cases}
\end{align*}
%***************************************************************************************************
Further define, $\displaystyle\mathcal{Q} \define \{\denom(c_i)\}_{i=1}^{n}\text{ where } c_i \in (\c_I,\c_C)$, $\mathcal{B}$ as the set of all possible bases of $\A_C$ such that $\A_\B$, $\B \in \mathcal{B}$ is invertible and $$\kappa = \lcm(\lcm(\mathcal{Q}),\,\abs{\lcm_{B\in\mathcal{B}}(\det(B))},\,\denom(\phi(\bz))),$$ where $\lcm$ stands for least common multiple. Observe that the quantities $\kappa^2\c_{\B_\u} \A_{\B_\u}^{-1}$, $\kappa^2 \phi(\bz)$, and $\kappa^2 \c_I$ are all integers. Therefore multiplying \cref{eq:12} with $\kappa^2$ we obtain
\begin{align*}
    \abs{\kappa^2 \c_I^\top\x^\u_I + \kappa^2 \c_{\B_\u}^\top\A_{\B_\u}^{-1}(\b-\A_I\x^\u_I) -\kappa^2 \phi(\bz) + \kappa^2 \c_{\B_\u}^\top\A_{\B_\u}^{-1}\u} <  \kappa^2\epsilon\label{eq:13}
\end{align*}
Let $\tilde\epsilon = \kappa^2\epsilon$. 
\begin{align}
    \tilde\epsilon & > \abs{\kappa^2 \c_I^\top\x^\u_I + \kappa^2 \c_{\B_\u}^\top\A_{\B_\u}^{-1}(\b-\A_I\x^\u_I) -\kappa^2 \phi(\bz) + \kappa^2 \c_{\B_\u}^\top\A_{\B_\u}^{-1}\u}\nonumber\\
    & \geq \abs{\abs{\kappa^2 \c_I^\top\x^\u_I + \kappa^2 \c_{\B_\u}^\top\A_{\B_\u}^{-1}(\b-\A_I\x^\u_I) -\kappa^2 \phi(\bz)} - \abs{\kappa^2 \c_{\B_\u}^\top\A_{\B_\u}^{-1}\u}}  \label{eq:14}
\end{align}
There are two terms in \cref{eq:14} only one of which involves $\u$, which can be upper bounded as,
\begin{align*}
    \abs{\kappa^2 \c_{\B_\u}^\top\A_{\B_\u}^{-1}\u} \leq \kappa^2 \norm{\c_C}\;\; \max_{B\in \mathcal{B}}\norm{\A^{-1}_{B}}\;\; \norm{\u}< \kappa^2  \max_{\mathcal{B}}\, \norm{\c_C}\,\norm{\A^{-1}_{B}} \delta.
\end{align*}
Letting $\beta:= \max_{B\in\mathcal{B}}\norm{\A^{-1}_B}$ and $K  := \max(\kappa^2\,\beta \norm{\c_C}, 1)$ yields, 
\begin{equation*}
    \abs{\kappa^2 \c_{\B_\u}^\top\A_{\B_\u}^{-1}\u} < K\delta
\end{equation*}
Now the remaining term from \cref{eq:14},
$\kappa^2 \c_I^\top\x^\u_I + \kappa^2 \c_{\B_\u}^\top\A_{\B_\u}^{-1}(\b-\A_I\x^\u_I) -\kappa^2 \phi(\bz)$, as highlighted earlier, is an integer. If this term is non-zero then 
\begin{equation}
    \abs{\kappa^2 \c_I^\top\x^\u_I + \kappa^2 \c_{\B_\u}^\top\A_{\B_\u}^{-1}(\b-\A_I\x^\u_I) -\kappa^2 \phi(\bz)} \geq 1 \ . \label{eq:16(i)}
\end{equation}
Letting $\tilde\delta = \min(\tilde\epsilon,\delta)$ and  $\tilde\epsilon < \frac{1}{2K}$, we have $\tilde\delta < \frac{1}{2K}$ implying that $K\tilde\delta  < 1/2.$ Now, we have
\begin{equation}
    \abs{\kappa^2\c_{\B_\u}\A_{\B_\u}\u}<K\tilde\delta < 1/2 \label{eq:16(ii)}
\end{equation}
Combining the inequalities \cref{eq:16(i)} and \cref{eq:16(ii)} yields,
\begin{align}
    \abs{\kappa^2 \c_I^\top\x^\u_I + \kappa^2 \c_{\B_\u}^\top\A_{\B_\u}^{-1}(\b-\A_I\x^\u_I) -\kappa^2 \phi(\bz)} - \abs{\kappa^2\c_{\B_\u}\A_{\B_\u}\u} > 1/2\nonumber\\
        \biggl\lvert{\abs{\kappa^2 \c_I^\top\x^\u_I + \kappa^2 \c_{\B_\u}^\top\A_{\B_\u}^{-1}(\b-\A_I\x^\u_I) -\kappa^2 \phi(\bz)} - \abs{\kappa^2\c_{\B_\u}\A_{\B_\u}\u}}\biggr\rvert > \frac{1}{2}\label{eq:17}
\end{align}
If $ \abs{\kappa^2 \c_I^\top\x^\u_I + \kappa^2 \c_{\B_\u}^\top\A_{\B_\u}^{-1}(\b-\A_I\x^\u_I)}$ is non-zero for any $\u \in \mathcal{N}_{\epsilon}(\bz)\cap \ocalU$ then from \cref{eq:14} and \cref{eq:17} one has
\begin{align*}
    1/2 & < \abs{\abs{\kappa^2 \c_I^\top\x^\u_I + \kappa^2 \c_{\B_\u}^\top\A_{\B_\u}^{-1}(\b-\A_I\x^\u_I) -\kappa^2 \phi(\bz)} - \abs{\kappa^2 \c_{\B_\u}^\top\A_{\B_\u}^{-1}\u}}\\
    & \leq \abs{\kappa^2 \c_I^\top\x^\u_I + \kappa^2 \c_{\B_\u}^\top\A_{\B_\u}^{-1}(\b-\A_I\x^\u_I) -\kappa^2 \phi(\bz) - \kappa^2 \c_{\B_\u}^\top\A_{\B_\u}^{-1}\u}\\
    & <  1/(2K)
\end{align*}
which presents a contradiction, since $K \geq 1$.

This implies that for every $u\in \mathcal{N}_{\delta}(\bz) \cap \ocalU$ where $0 < \delta < 1/(2K)$ we have 
$$\kappa^2 \c_I^\top\x^\u_I + \kappa^2 \c_{\B_\u}^\top\A_{\B_\u}^{-1}(\b-\A_I\x^\u_I) -\kappa^2 \phi(\bz) = 0.$$ It follows from \cref{eq:11} that $\exists\, \delta > 0$ such that $\forall $ $\u\in \mathcal{N}_{\delta}(\bz)\cap \overline{U}$ we have $\phi(\bz)-\phi(\u) = \c_{\B_\u}^\top\A_{\B_\u}^{-1}\u$. It follows,
\begin{align*}
    \abs{\phi(\u)-\phi(\bz)} = \abs{\c_{\B_\u}^\top\A_{\B_\u}^{-1}\u}
    \leq \norm{\c_{\B_\u}}\cdot\norm{\A^{-1}_{\B_\u}}\cdot\norm{\u}
    \leq \beta \norm{\c_C}\, \norm{\u}
\end{align*}
As a result,
\begin{align*}
    \abs{\phi(\u)-\phi(\bz)} \leq \Gamma \norm{\u}\;\;
    \text{where}\;\; \Gamma = \beta \norm{\c_C}\;\;
\end{align*}
\end{proof}
\begin{corollary}
Consider the value function $\phi(\u)$ of a mixed integer linear programming problem,
\begin{optprog*}
    {$\phi(\u) = $ minimize}  & \objective{\hspace{5mm}\c_I^{\top}\x_I + \c_C^{\top}\x_C}\\
    \hspace{2mm} subject to & \A_I\,\x_I + \A_C\x_C = \b + \u\\
     & \E_I\x_I + \E_C\x_C \leq \f\hspace{7mm}\\
     & \x_I\in \Z^{n_1}\hspace{1.2cm}\\ 
     & \x_C \in \R^{n_2}\hspace{1.2cm}
\end{optprog*}
There exists a $\delta > 0$ such that for every $ \u\in \mathcal{N}_{\delta}(\bz)\cap \calU$, we have $\phi(\bz) \leq \phi(\u) + \Gamma\norm{\u}$
where $\Gamma$ is a constant which depends on $\A_C$ and $\c_C$ and $\c_I$ .
\end{corollary} 
\begin{proof}
Without loss of generality, we can represent the given mixed integer linear program in the following form,
\begin{optprog*}
    {$\phi(\u') = $ minimize}  & \objective{\hspace{5mm}{\c'_I}^{\top}\x'_I + {\c'_C}^{\top}\x'_C}\\
    \hspace{2mm} subject to & \A'_I\,\x'_I + \A'_C\x'_C = \b' + \u'\\
     & \x'_I \in \Z_{+}^{n_1}\hspace{1.5cm}\\ 
     & \x'_C \in \R_{+}^{n_2}\hspace{1.5cm}
\end{optprog*}
where ${\c'_I}^\top = (\c_I, -\c_I)$, ${\c'_C}^\top = (\c_C, -\c_C,\bz)$ and
\begin{align*}
 \A'_I =    \begin{pmatrix}
\A_I & -\A_I\\
\E_I & -\E_I
\end{pmatrix},
&
\A'_C = \begin{pmatrix}
\A_I & -\A_C &\mathbb{O}\\
\E_I &  -\E_I&\mathbb{I}
\end{pmatrix},
\x'_I = \begin{pmatrix}
\x^{+}_I\\
\x^{-}_I
\end{pmatrix},
\x'_C = \begin{pmatrix}
\x^{+}_C\\
\x^{-}_C\\
\s
\end{pmatrix},
\\
&
\b' = \begin{pmatrix}
\b\\
\f
\end{pmatrix},
\u' = \begin{pmatrix}
\u\\
\bz
\end{pmatrix}
\end{align*}
and $\mathbb{O}$, $\mathbb{I}$ denote the matrix of all zeros and identity matrix of appropriate dimensions. The result follows from \cref{prop:4}.
\end{proof}
\cref{thm:1} follows consequently.
\subsection{Proof of \cref{thm:2}}
\begin{proposition}\label{prop:5}
Consider the value function $\phi(\u)$ of a mixed integer quadratic programming problem,
\begin{optprog*}
    {$\phi(\u) = $ minimize}  & \objective{\frac{1}{2}\x_C^{\top}Q_{CC}\x_C + \frac{1}{2}\x_I^{\top}Q_{II}\x_I + \x_I^{\top}Q_{IC}\x_C- {\c}_I^\top\x_I - {\c}_C^\top\x_C}\\
    subject to & \hspace{1.2cm} \A_I\,\x_I + \A_C\x_C = \b + \u\\
    & \norm{\x_I}_{\infty}\leq M \hspace{0.5cm}\\
     & \x_I\in \Z^{n_1}\hspace{0.4cm}\\ 
     & \x_C \in \R_{+}^{n_2}\hspace{0.4cm}
\end{optprog*}
There exists a $\delta > 0$ such that for every $ \u\in \mathcal{N}_{\delta}(\bz)\cap \calU$, we have $\phi(\bz) \leq \phi(\u) + K_1\norm{u}^2 + K_2\norm{u}$ for some $K_1, K_2$ which depend on $\A,\Q,{\c}$ and $M$.
\end{proposition}
\begin{proof}
 As the set $S_I = \set{\x_I \in \Z^{n_1}}{\norm{\x_I}_{\infty} \leq M}$ is finite, let $S_I = \{\x^{(i)}_I\}_{i \in [k]}$ for some natural number $k$. We now characterize the value functions of the continuous restrictions, parameterized in $\x_I^{(i)},\,i \in [k]$ as,
 \begin{optprog*}
    {$\Phi(\u,\x^{(i)}_I) = $ minimize}  & \objective{\frac{1}{2}{\x^{(i)}_I}^{\top}Q_{II}\x^{(i)}_I + \frac{1}{2}{\x_C}^{\top}Q_{CC}{\x_C} + {\x^{(i)}_{I}}^{\top}Q_{CI}{\x_C} -{\c}_I^\top\x^{(i)}_{I}- {\c}_C^\top{\x_C}}\\
    $\QPui$ \hspace{2mm}subject to & \hspace{1.6cm} \A_C\,\x_C = \b + \u - \A_I\x^{(i)}_{I}\\
     & \x_C \in \R_{+}^{n_2}\hspace{1.8cm}
\end{optprog*}
and $\phi(\u) = \min_{1 \leq i \leq k} \,\Phi(\u,\x^{(i)}_I) $. We can assume without loss of generality that $\QPui$ is feasible for all $((\u,\x^{(i)}_I)) \in \calU\times S_I$, since, if $\exists$ $((\u,\x^{(i)}_I)) \in \calU\times S_I$, such that $\QPui$ is infeasible, we can assign $\Phi(\u,\x^{(i)}_I) = \infty$. Observe that continuous quadratic programs are lower semi-continuous at $\u = \bz$ (\cite{bertsekas2003convex}, Proposition 6.5.2) and the minimum of a finite number of lower semi-continuous functions is lower semi-continuous as well. This implies that $\phi(\u)$ is lower semi-continuous at $\u = \bz$ and consequently $\lim_{\u\to 0} \phi(\u) = \phi(\bz)$ where $\u\in \ocalU$. The first order necessary conditions for optimality (KKT conditions) for $\QPui$ can be expressed as,
\begin{align}\label{eq:KKT}
\begin{split}
    \begin{bmatrix}
        Q_{CC} & -\A_C^\top & \A_C^\top & -I\\ 
        \A_{C}& \bz & \bz & \bz
    \end{bmatrix} 
    \begin{bmatrix}
        \x_C \\
        \blambda_{+} \\
        \blambda_{-}\\
        \bpi
    \end{bmatrix} 
    & = 
    \begin{bmatrix}
        {\c}_C - Q^\top_{CI}\x^{(i)}_I\\
        \b + \u - \A_I\x^{(i)}_I
    \end{bmatrix}\\
    \boldsymbol{\pi}^{\top}\x_C = 0\\
    \x_C,\blambda_{+},\blambda_{-}, \bpi & \geq \bz
\end{split}
\end{align}

If $(\x_C,\x^{(i)}_I)$ is a solution to the KKT conditions then it satisfies the complementary slackness conditions. Let $J=[n_2]$ and consider the partition of $J$ for $i \in [k]$, $J(i) = \{J^{(i)}_>,J^{(i)}_=\}$ such that $J^{(i)}_{=} = \set{j\in J}{\boldsymbol{\pi}^{(j)} = 0 }$ and $J^{(i)}_> = J\backslash J^{(i)}_{=}$. It follows that $\forall j \in J_>^{(i)}$, $x_C^{(j)} = 0$. Define,
\begin{align*}
    \A_{aug} =     \begin{bmatrix}
        Q_{CC} & -\A_C^\top & \A_C^\top & -I\\ 
        \A_{C}& \bz & \bz & \bz
    \end{bmatrix} 
\end{align*}
and let the columns of $\A_{aug}$ be indexed by $ [2n_2+2m]$. Define $\A_{aug}^{J(i)}$ as a sub-matrix $\A_{aug}$ which has all the columns of $\A_{aug}$ except the columns $ j \in J^{(i)}_> \cup ((n_2+2m)+J^{(i)}_=).$

Consider the set of solutions to \cref{eq:KKT},
$$
    P_{J(i)} \define \Bigg\{(\x_C,\blambda_{+},\blambda_{-},\bpi) \in \R^{n_2}_+\times\R^{m}_+\times\R^{m}_+\times\R^{n_2}_+ :\hspace*{4cm}
$$
$$
    \hspace*{4cm}\left.\A^{J(i)}_{aug}
    \begin{bmatrix}
        \x_C \\
        \blambda_{+} \\
        \blambda_{-}\\
        \bpi
    \end{bmatrix} = \begin{bmatrix}
        {\c}_C - Q^\top_{CI}\x^{(i)}_I\\
        \b + \u - \A_I\x^{(i)}_I
    \end{bmatrix}
    \right\}
$$
Observe that as the polyhedron $P_{J(i)}$ doesn't contain a line, it must have an extreme point. Consider an extreme point of $P_J$ corresponding to the basis $\B_\u^{(i)}$ of columns of $\A^{J(i)}_{aug}$. Let $\A_{\B_\u}^{(i)}$ be the sub-matrix formed by the columns corresponding to basis $\B_\u^{(i)}$.\\

In any solution to \cref{eq:KKT} we have,
\begin{equation}
    \x_C = {\C}^{\u}_1 ({\c}_C - Q^\top_{CI}\x^{(i)}_I) + {\C}^{\u}_2(\b + \u - \A_I\x^{(i)}_I)\label{eq:23}
\end{equation}
where ${\C}^{\u}_1$ and ${\C}^{\u}_2$ are submatrices of $\A_{\B_\u}^{(i){-1}}$. Substituting $\x_C$ from \cref{eq:23} in the objective function of $\QPui$, we obtain
\begin{equation}
    \Phi(\u,\x^{(i)}_I) = \c_{\u1}^\top\u  + \frac{1}{2}\u^\top\C_2^\u Q_{CC}\C_2^u\u + \Theta(\b,\x^{(i)}_I,\C_1^{\u},\C_2^{\u}) \label{def:Phi}
\end{equation}
where,
\begin{align*}
    \c_{\u1} & = 2{\zeta^{\u}}^\top Q_{CC}\C^{\u}_2 + \x^{(i)\top}_IQ_{CI}\C^{\u}_2 + \c^\top_C\C^{\u}_2,\\
    \Theta(\b,\x^{(i)}_I,\C_1^{\u},\C_2^{\u}) & = \frac{1}{2}\x^{(i)\top}_IQ_{II}\x^{(i)}_I + \c_C^\top{\zeta^{\u}} + \frac{1}{2}{\zeta^{\u}}^\top Q_{CC} {\zeta^{\u}} + {\x^{(i)\top}_I} Q_{CI}\zeta^{\u},  \text{ and}\\
    \zeta^{\u} & = \left(\C^{\u}_1 \left(\c_C - Q^\top_{CI}\x^{(i)}_I\right) + \C^{\u}_2\left(\b - \A_I\x^{(i)}_I\right)\right)
\end{align*}
Now, for every $\epsilon > 0$ there exists a $\delta > 0$ such that
\begin{equation*}
    % \abs{\phi(\u) - \phi(\bz)} < \epsilon \hspace{1cm} \forall \norm{\u} < \delta.
        \abs{\phi(\u) - \phi(\bz)} < \epsilon \hspace{1cm} \forall \mathcal{N}_{\delta}(\mathbf{0})\cap \overline{U}
\end{equation*}
It follows that $\forall \norm{\u} < \delta$, we have
\begin{align}
    \abs{\c_{\u1}^\top\u  + \frac{1}{2}\u^\top\C_2^\u Q_{CC}\C_2^\u\u + \Theta(\b,\x^{(i)}_I,\C_1^{\u},\C_2^{\u}) - \phi(\bz)} &< \epsilon, \text{ i.e.}\nonumber\\
    \abs{\abs{\Theta(\b,\x^{(i)}_I,\C_1^{\u},\C_2^{\u}) - \phi(\bz)}- \abs{\c_{\u1}^\top\u  + \frac{1}{2}\u^\top\C_2^\u Q_{CC}\C_2^\u\u}} &< \epsilon \label{eq:27}
\end{align}

Let $\A_\B^* = \arg\max\set{\norm{\A^{-1}_\B}}{\B\in \mathcal{B}}$ where $\mathcal{B}$ is the set of all possible bases of $\A_{aug}$.

Since $\norm{\A_{\B_\u}^{(i){-1}}}_F \leq \norm{{\A_\B^*}^{-1}}_F$ and $\C^{\u}_1$ and $\C^{\u}_2$ are submatrices of $\A_{\B_\u}^{(i){-1}}$ it follows that $\norm{\C^{\u}_2}_F \leq \norm{\A_{\B_\u}^{(i){-1}}}_F$ and $\norm{\C^{\u}_1}_F \leq \norm{\A_{\B_\u}^{(i){-1}}}_F$. We can bound the second term in inequality \cref{eq:27} as\\
\begin{align}
    \abs{\c_{\u1}^\top\u  + \frac{1}{2}\u^\top\C_2^\u Q_{CC}\C_2^\u\u} \leq K_1\norm{\u}^2 + K_2\norm{\u} \label{ineq:main}
\end{align}
where $K_1 = \norm{Q}_F\norm{{\A_\B^*}^{-1}}_F^2$ and 
\begin{align*}
    \norm{\c_{\u1}} & = \norm{2{\zeta^{\u}}^{\top}Q_{CC}\C^{\u}_2 + \x^{(i)\top}_IQ_{CI}\C^{\u}_2 + \c^{\top}_C\C^{\u}_2}\\
    & \leq 2\left(\norm{\C^{\u}_1}_F\norm{\c_C - Q^{\top}_{CI}\x^{(i)}_I} + \norm{\C^{\u}_2}_F\,\norm{\b - \A_I\x^{(i)}_I}\right)\left(\norm{Q_{CC}}_F\norm{\C^{\u}_2}_F\right)\\
    & \hspace{1cm} + \norm{\x^{(i)}_I}\,\norm{Q_{CI}}_F\,\norm{\C^{\u}_2}_F + \norm{\c_C}\,\norm{\C^{\u}_2}_F\\
    & \leq 2\left(\norm{\c} + M\norm{Q}_F + \norm{\b} + M \norm{\A_{I}}_F\right)\norm{Q}_F\,\norm{{\A_\B^*}^{-1}}^2_F\\
    & \hspace{1cm} + M\norm{Q}_F \norm{{\A_\B^*}^{-1}}_F + \norm{\c}\, \norm{{\A_\B^*}^{-1}}_F = K_2.
\end{align*}
Since $\norm{u} < \delta$ it follows that $\abs{\c_{\u1}^\top\u  + \frac{1}{2}\u^\top\C_2^\u Q_{CC}\C_2^\u\u} < K_1\delta^2 + K_2\delta$. Define, 
$$\psi \define \min_{\x^{(i)}, \C^{\u}_1, \C^{\u}_2} \abs{\Theta(\b,\x^{(i)}_I,\C_1^{\u},\C_2^{\u}) - \phi(\bz)}.
$$
Since $\x^{(i)}_I$ is finite, $\b$ is fixed and $\abs{\mathcal{B}}$ is finite, $\psi$ exists. Assume $\psi\neq 0$. Without loss of generality, take $\delta \leq \epsilon$ and take $\epsilon <\psi/2$ such that $K_1\delta^2 + K_2 \delta < \dfrac{\psi}{2}$. Now, we have
\begin{align}
    \abs{\Theta(\b,\x^{(i)}_I,\C_1^{\u},\C_2^{\u}) - \phi(\bz)} \geq \psi \text{ and } \abs{\c_{\u1}^\top\u  + \frac{1}{2}\u^\top\C_2^\u Q_{CC}\C_2^\u\u} < \frac{\psi}{2}\nonumber\\
   \abs{\Theta(\b,\x^{(i)}_I,\C_1^{\u},\C_2^{\u}) - \phi(\bz)} - \abs{\c_{\u1}^\top\u  + \frac{1}{2}\u^\top\C_2^\u Q_{CC}\C_2^\u\u} > \frac{\psi}{2}\nonumber\\
   \abs{\abs{\Theta(\b,\x^{(i)}_I,\C_1^{\u},\C_2^{\u}) - \phi(\bz)} -  \abs{\c_{\u1}^\top\u  + \frac{1}{2}\u^\top\C_2^\u Q_{CC}\C_2^\u\u}} > \frac{\psi}{2} \label{eq:28}
\end{align}
From \cref{eq:27} and \cref{eq:28} we get
\begin{align*}
     \psi/2 & < \abs{\abs{\Theta(\b,\x^{(i)}_I,\C_1^{\u},\C_2^{\u}) - \phi(\bz)} -  \abs{\c_{\u1}^\top\u  + \frac{1}{2}\u^\top\C_2^\u Q_{CC}\C_2^\u\u}} < \epsilon < \psi/2
\end{align*}
This yields a contradiction. Hence $\psi = 0$, which further yields $\Theta(\b,\x^{(i)}_I,\C_1^{\u},\C_2^{\u}) = \phi(\bz)$. Substituting in \cref{def:Phi}, it follows from \cref{ineq:main} that
\begin{align*}
    % \abs{\phi(\u)-\phi(\bz)} < K_1\norm{u}^2 + K_2\norm{u} \label{eq:29}
   \phi(\bz) \leq \phi(\u)+ K_1\norm{\u}^2 + K_2\norm{\u} \label{eq:29}
\end{align*}

\end{proof}
\begin{corollary}\label{cor:6}
Consider the value function $\phi(\u)$ of a mixed integer quadratic programming problem,
\begin{optprog*}
    {$\phi(\u) = $ minimize}  & \objective{\frac{1}{2}\x_C^{\top}Q_{CC}\x_C + \frac{1}{2}\x_I^{\top}Q_{II}\x_I + \x_I^{\top}Q_{IC}\x_C- {\c}_I^\top\x_I - {\c}_C^\top\x_C}\\
    subject to & \hspace{1.2cm} \A_I\,\x_I + \A_C\x_C = \b + \u\\
    & \E_I\x_I + \E_C\x_C \leq \f \hspace{7mm}\\
    & \norm{\x_I}_{\infty}\leq M \hspace{0.5cm}\\
     & \x_I\in \Z^{n_1}\hspace{0.4cm}\\ 
     & \x_C \in \R^{n_2}\hspace{0.4cm}
\end{optprog*}
There exists a $\delta > 0$ such that $\forall$ $ \u\in \mathcal{N}_{\delta}(\bz)\cap \calU$, $\phi(\bz) \leq \phi(\u) +  K_1\norm{\u}^2 + K_2\norm{\u}$ for some $K_1, K_2$ which depend on $\A,\Q,{\c}$ and $M$.
\end{corollary}
\begin{proof}
Without loss of generality, we can represent the given mixed integer quadratic program in the following form,
\begin{optprog*}
    minimize & \objective{\frac{1}{2}{\x'_C}^{\top}{Q'}_{CC}\x'_C + \frac{1}{2}\x_I^{\top}Q_{II}\x_I + \x_I^{\top}{Q'}_{IC}\x'_C- {\c}_I^\top\x_I - {\c'}_C^\top\x'_C}\\
    subject to & \hspace{1.2cm} \A_I\,\x_I + {\A'}_C\x'_C = {\b'} + {\u'}\\
    & \norm{\x_I}_{\infty}\leq M \hspace{0.5cm}\\
     & \x_I\in \Z^{n_1}\hspace{0.4cm}\\ 
     & \x'_C \in \R_{+}^{n_2}\hspace{0.4cm}
\end{optprog*}
where
\begin{align*}
& Q'_{CC} = \begin{pmatrix}
Q_{CC} & -Q_{CC} & \bz\\
-Q_{CC} & Q_{CC} & \bz\\
\bz & \bz & \bz
\end{pmatrix},
{\Q'}_{IC} = \begin{pmatrix}
   \Q_{IC} & -\Q_{IC} & \bz
\end{pmatrix},\\
& \A'_C = \begin{pmatrix}
\A_I & -\A_C &\mathbb{O}\\
\E_I &  -\E_I&\mathbb{I}
\end{pmatrix},
\x'_C = \begin{pmatrix}
\x^{+}_C\\
\x^{-}_C\\
\s
\end{pmatrix},
\c'_C = \begin{pmatrix}
\c_C\\ 
-\c_C\\
\bz
\end{pmatrix},
\b' = \begin{pmatrix}
\b\\
\f
\end{pmatrix},
\u' = \begin{pmatrix}
\u\\
\bz
\end{pmatrix}
\end{align*}
and $\mathbb{O}$, $\mathbb{I}$ denote the matrix of all zeros and identity matrix of appropriate dimensions. The result follows from \cref{prop:5}.
\end{proof}
\cref{thm:2} follows immediately.
%\section{Exact Penalty Representation for Mixed Integer Convex Programming Problems}
\subsection{Proof of \cref{thm:3}}
Consider the following mixed integer convex programming problem.
\begin{optprog*}
    {$\hat{\phi}(\u) = $ minimize}  & \objective{\hspace{1.3cm}f({\x_I},{\x_C})}\\
    $\textup{\MICP}_\u$ \qquad subject to & \hspace{0cm} \A_I\,\x_I + \A_C\x_C = \b + \u\\
    & \E_I\x_I + \E_C\x_C \leq \f \hspace{7mm}\\
    & \norm{\x_I}_{\infty}\leq M \hspace{0.5cm}\\
    % & \x_I\in \Z_{+}^{n_1}\hspace{0.4cm}\\ 
    & (\x_I, \x_C) \in \Z^{n_1}\times\R^{n_2}\hspace{0.3cm}
\end{optprog*}
where $0 \leq M < \infty$. Denote by $F_{\u} \define X \cap \Hu \cap S_I $ the feasible region of $\textup{\MICP}_\u$, where the set $S_I = \set{\x_I \in \Z^{n_1}}{\norm{\x_I}_{\infty} \leq M}$ is finite.  As before, we let $S_I = \{\x^{(1)}_I, \x_I^{(2)}, \ldots, \x_I^{(k)}\}$. Using this enumeration we characterize the value functions of the continuous restrictions, parameterized in $\x_I^{(i)},\,i \in [k]$ as,
\begin{optprog*}
    {$\Phi(\u,\x^{(i)}_I) = $ minimize}  & \objective{\hspace{6mm}f(\x^{(i)}_I,\x_C)}\\
    $\CPui$ \qquad subject to & \A_C\,\x_C = \b + \u - \A_I\x^{(i)}_{I}\\
    & \E_C\x_C \leq \f - \E_I\x^{(i)}_{I}\hspace{8mm}\\
     & \x_C \in \R^{n_2}\hspace{1.8cm}
\end{optprog*}
where $\Phi(\u,\x^{(i)}_I)$ is the value function of $\CPui$, a continuous convex  optimization problem. Since the constraints of $\CPui$ are linear, Slater's condition holds and therefore the value function is lower semi-continuous.
\begin{lemma}\label{lem:6}
If the integer variables are bounded, then the value function $\displaystyle \hat{\phi}(\u) = \min_{i \in [k]} \Phi(\u,\x^{(i)}_I)$ is lower semi-continuous at $\bz$.
\end{lemma}
\begin{proof}
Observe that for $\u \in \calU, \x^{(i)}_I \in S_I$, $\Phi(\u,\x^{(i)}_I)$ is lower semi-continuous at ${\u}$ if and only if strong duality holds at ${\u}$ [Proposition 6.5.2,\cite{bertsekas2003convex}]. If for some $\x^{(i)}_I$, $\Phi(\bz, \x^{(i)}_I) < +\infty$ then $\Phi(\u, \x^{(i)}_I)$ is lower semi-continuous at $\u = \bz$ due to Slater Condition (in case of affine constraints only feasibility is required). On the contrary, $\Phi(\bz,\x^{(i)}_I) = + \infty$ implies that the polyhedron $P = \set{\x \in \R^{n_2}}{\tilde{\A}\x = \tilde{\b}, \,\,\x \geq \bz}$ is empty, where
\begin{align*}
    \tilde{\A} = \begin{pmatrix}
      {\A}_C & -{\A}_C & \bz\\
      {\E}_C & -{\E}_C & \mathbb{I}
    \end{pmatrix},~
    \x = \begin{pmatrix}
    \x_C^+\\
    \x_C^-\\
    \s
    \end{pmatrix},~
    \tilde{\b} = \begin{pmatrix}
     {\b}  -{\A}_I{\x}^{(i)}_I\\
     {\f} - {\E}_I{\x}^{(i)}_I
    \end{pmatrix}.
\end{align*}
Let $C$ be the closed convex cone spanned by the columns of $\tilde{\A}$. If $P = \emptyset$ then $\tilde{\b} \not\in C$. Since $C$ is a closed, 
$\exists\,\,\delta_{i} > 0$ such that $\mathcal{N}_{\delta_i}(\tilde{\b}) \cap C = \emptyset$. Consider $\u \in \R^m \cap \mathcal{N}_{\delta_i}(\bz)$, then $\tilde{\b} + \begin{pmatrix}\u\\\bz\end{pmatrix} \not \in C$.

For $J = \set{i\in [k]}{\Phi({\x}^{(i)}_I,\bz) < \infty}$, define $\displaystyle\delta_{\min} \define \min_{i \in [k]\setminus J}\delta_i$.

Consider y such that $y  < \hat{\phi}(\bz)$, it follows that 
\begin{align*}
    y & < \hat{\phi}(\bz) < \min_{i \in J}\Phi(\bz,{\x}^{(i)}_I) \leq \Phi(\bz,{\x}^{(i)}_I) .
\end{align*}
However, since every $\Phi({\u},{\x}^{(i)}_I), \, i \in J$ is lower semi-continuous at $\bz$, there exists a $\delta_i > 0$ such that $\Phi({\u},{\x}^{(i)}_I) > y$ for every ${\u} \in \mathcal{N}_{\delta_i}(\bz), \, i \in J$. Let $\delta'_{\min} = \min(\delta_i),\, i \in J$ and $\delta'' = \min(\delta_{\min},\delta'_{\min})$. This means that for every ${\u} \in \mathcal{N}_{\delta''}(\bz)$ $\Phi({\u},{\x}^{(i)}_I) > y \implies  \hat{\phi}({\u}) > y$. The result follows.
\end{proof}
Consider $\ocalU$ as defined in \cref{cor:4}, it follows from \cref{lem:6} that $\lim_{\u\to \bz} \hat{\phi}(\u) = \hat{\phi}(\bz)$ where $\u \in \ocalU$. Now, for every $\epsilon > 0$ there exists a $\delta > 0$ such that
\begin{equation*}
   \abs{\hat{\phi}(\u) - \hat{\phi}(\bz)} < \epsilon\,\, \forall \u\,\, \in \mathcal{N}_{\delta}(\bz) \cap \ocalU
\end{equation*}
\begin{lemma}\label{lem:8}
Consider the value function $\hat{\phi}(\u)$ as defined earlier for $\textup{\MICP}_\u$. There exists a $\delta > 0$ such that for all $\u \in \mathcal{N}_\delta(\bz) \cap \ocalU$, $\displaystyle \hat{\phi}(\u) = \min_{\x_I\in S_{=}^{\bz}} \Phi(\u,\x_I)$ where $S_{=}^{\bz} = \set{\x_I \in S_I}{\Phi(\bz,\x_I) = \hat{\phi}(\bz)}$
\end{lemma}
\begin{proof}
We can partition the index sets corresponding to solutions of $\CPui$ into two sets, 
$$
S_{>}^{\bz} = \set{\x_I\in S_I}{\Phi(\bz,\x_I) > \hat{\phi}(\bz)} \text{ and } S_{=}^{\bz} = \set{\x_I\in S_I}{\Phi(\bz,\x_I) = \hat{\phi}(\bz)}.
$$
For each $\x^{(i)}_I \in S_I$, define $ \overline{\calU}_i=\set{\u \in \calU}{\Phi({\u},\x^{(i)}_I)\leq \Phi(\bz,\x^{(i)}_I)} $. If there is no limiting sequence to ${\bz}$ in $\overline{\calU}_i$ then there exists a neighbourhood $\delta_i$ such that $\mathcal{N}_{\delta_i}(\bz) \cap \overline{\calU}_i = \emptyset$. Alternatively, if there does exist a limiting sequence to $\bz$ and since $\Phi({\u},{\x}_I^{(i)})$ is lower semi-continuous,
\begin{equation*}
    \lim_{{\u}\to \bz} \Phi({\u},{\x}^{(i)}_I) = \Phi(\bz,{\x}^{(i)}_I)
\end{equation*}
Thus, for every $\epsilon_i > 0$ there exists a $\delta_i > 0$ such that
\begin{align*}
    \abs{\Phi({\u},{\x}^{(i)}_I) - \Phi(\bz,{\x}^{(i)}_I)}< \epsilon_i \,\,\,\,\, \forall \,\,\,\,{\u}\in \overline{\calU}_i\cap \mathcal{N}_{\delta_i}(\bz)
\end{align*}
For $\x^{(i)}_I \in S_{>}^{\bz}$, let $\epsilon_i = \Phi(\bz,{\x}^{(i)}_I) - \hat{\phi}(\bz)$. Assuming that $\Phi(\bz,{\x}_I^{(i)})$ is finite, $\forall$ $\u \in \overline{\calU}_i\cap \mathcal{N}_{\delta_i}(\bz)$
\begin{align*}
    \Phi(\bz,{\x}^{(i)}_I) - \Phi({\u},{\x}_I^{(i)}) & <  \Phi(\bz,{\x}^{(i)}_I) - \hat{\phi}(\bz)\\
    \hat{\phi}(\bz) & < \Phi({\u},{\x}^{(i)}_I)
\end{align*}
Thus, for all $\x^{(i)}_I\in S_{>}^{\bz}$ there exists a $\delta_i$ such that $\Phi({\u},{\x}_I^{(i)}) > \hat{\phi}(\bz)$. Assigning $\delta_{\min} = \min_{\x^{(i)}_I\in S^{\bz}_{>}} \delta_i$ yields $\displaystyle \hat{\phi}({\u}) = \min_{\x_I\in S^{\bz}_{=}}\Phi({\u},{\x}_I)$ for every $\u \in \mathcal{N}_{\delta_{\min}}(\bz)\cap \overline{\calU}$.
\end{proof}
\begin{lemma}\label{lem:9}
Consider the value function $\hat{\phi}(\u)$ as defined earlier for $\textup{\MICP}_\u$. There exists a $\delta > 0$ and $0 < \Gamma< \infty$ such that for every $ \u\in \mathcal{N}_{\delta}(\bz)\cap \calU$, $$\hat{\phi}(\bz)\leq \hat{\phi}(\u) + \Gamma\norm{\u}.$$
\end{lemma}
\begin{proof}
Consider $S^{\bz}_{=}$ as defined in \cref{lem:8}. Analogously, define $$S_{=}^{\u} \define \set{\x_I\in S_{=}^{\bz}}{\Phi(\u,\x_{I}) = \hat{\phi}(\u)}.$$ For $\x^{(i)}_I \in S^{\u}_{=}$, the first order necessary conditions for optimality (KKT conditions) for $\CPui$ can be expressed as, %For  writing the KKT conditions\\
\begin{align}\label{eq:32}
\begin{split}
    \nabla_{C}{f(\x_I^{(i)},\x_C}) & = \A^\top_C\blambda^{\u(i)}_{\A_C} - \E_C^\top\blambda^{\u(i)}_{\E_C}\\
    \A_C\x_C & = \b + \u - \A_I\x^{(i)}_I\\
    \E_C\x_C+\E_I\x^{(i)}_I & \leq \f\\
    {\blambda^{\u(i)}_{\E_C}}^\top(\f - \E_C\x_C - \E_I\x^{(i)}_I) & = \bz\\
    {\blambda^{\u(i)}_{\A_C}}^\top(\A_C\x_C + \A_I\x^{(i)}_I - \b - \u) & = \bz\\
    {\blambda^{\u(i)}_{\E_C}} & \geq \bz
\end{split}
\end{align}
Define $H^{(j)}_{\u} = \set{\x_C \in \R^{n_2}}{\A_C\x_C = \b +\u- \A_I\x^{(j)}_I}$. Let $(\x_I^{(j)},\x_C) \in X \cap H^{(j)}_{\u}$ where $\x_I^{(j)} \in S_{=}^{\u}$.
Additionally, let the corresponding Lagrange multipliers at $\u = \bz$ for $\x^{(j)}_I$ be $\blambda^{(j)}_{\A_C}$ and $\blambda^{(j)}_{\E_C}$.
\begin{align*}
    f(\x^{(j)}_I,\x^{\u}_C) & \geq f(\x^{(j)}_I,\x^{\bz^*}_C) +  \nabla_C{f(\x_I^{(j)}}\x^{\bz*}_C)^{\top}(\x^{\u}_C -  \x^{\bz*}_C )\\ 
    & \hspace{3cm}+\nabla_I{f(\x_I^{(j)}}, \x^{\bz*}_C)^{\top}(\x^{(j)}_I -  \x^{(j)}_I)\\
    & = f(\x^{(j)}_I,\x^{\bz^*}_C) + ( \blambda^{(j)^\top}_{\A_C}\A_C - \blambda^{(j)^\top}_{\E_C}\E_C) (\x^{\u}_C -  \x^{\bz*}_C)\\
    & = f(\x^{(j)}_I,\x^{\bz^*}_C) +  \blambda^{(j)^\top}_{\A_C}\A_C(\x^{\u}_C -  \x^{\bz*}_C) - \blambda^{(j)^\top}_{\E_C}\E_C (\x^{\u}_C -  \x^{\bz*}_C)\\
    & = f(\x^{(j)}_I,\x^{\bz^*}_C) +  \blambda^{(j)^\top}_{\A_C}(\A_C \x^{\u}_C -  \A_C \x^{\bz*}_C) - \blambda^{(j)^\top}_{\E_C}(\E_C \x^{\u}_C -  \E_C \x^{\bz*}_C)\\
    & = f(\x^{(j)}_I,\x^{\bz^*}_C) +  \blambda^{(j)^\top}_{\A_C}(\b+\u-\A_I\x^{(j)}_I - \b + \A_I\x^{(j)}_I)\\ 
    & \hspace{3cm}-\blambda^{(j)^\top}_{\E_C}(\E_C \x^{\u}_C -  \f + \E_I\x^{(j)}_I)
\end{align*}
where the equalities follow from the set of equations \cref{eq:32}. Furthermore, observe that since $\blambda^{(j)^\top}_{\E_C} \geq \bz$ and $\E_C\x^{\u}_C + \E\x^{(j)}_I \leq \f$, it follows that 
\begin{align*}
   f(\x^{(j)}_I,\x^{\u}_C) & \geq f(\x^{(j)}_I,\x^{\bz^*}_C) +  \blambda^{(j)^\top}_{\A_C}(\u), \\
      \inf_{X\cap H^{(j)}_{\u}}f(\x^{(j)}_I,\x^{\u}_C) & \geq f(\x^{(j)}_I,\x^{\bz^*}_C) +  \blambda^{(j)^\top}_{\A_C}(\u), \text{ i.e.,}\\
   \hat{\phi}({\u)} & \geq \hat{\phi}(\bz) +  \blambda^{(j)^\top}_{\A_C}(\u)
\end{align*}
Rearranging the terms in the above inequality,     $-\blambda^{(j)^\top}_{\A_C}\u \geq \hat{\phi}(\bz) - \hat{\phi}({\u)}$. As $\hat{\phi}(\bz)\geq \hat{\phi}(\u)$ $\forall\,\,\, \u \in \ocalU$, using Cauchy-Schwarz on the left-hand side of the inequality,
\begin{align*}
    \norm{\blambda^{(j)}_{\A_C}}\cdot\norm{\u} \geq \hat{\phi}(\bz) - \hat{\phi}({\u)} .
\end{align*}
Defining $\Gamma = \max_{i\in S^{\bz}_{=} }{\norm{\blambda^{(j)}_{\A_C}}}$ yields
$$
    \hat{\phi}(\bz) - \hat{\phi}(\u)\leq \Gamma\norm{\u}\,\,\forall\,\, \u \in \mathcal{N}_{\delta}(\bz)\cap\ocalU
$$
where $\delta = \delta_{min}$ as defined in the proof of \cref{lem:8}. The result follows.
\end{proof}
This completes the proof of \Cref{thm:3}(a). In the following section we present some cases where an equivalence can be established between (MICP) and (MICP) with bounded integer variables. In particular, we highlight that  if the objective function $f$ satisfies $\rec(f)\cap\rec(F_R)\setminus\{\bz\} = \emptyset$, then there exists $M < \infty$ such that the following equivalence holds.
$$ \min\set{f(\x)}{\A\x = \b, \x \in X} = \min\set{f(\x)}{\A\x = \b, \norm{\x}_{\infty} \leq M, \x \in X}.$$ Observe that if $f$ is $\mu$-strongly convex, then $\rec(f) = \emptyset$, as the level sets of $f$ are bounded. 

\subsubsection{MICPs with implicit integer boundedness}
\begin{proposition}\label{prop:13.1}
   If \textup{\MICP} is feasible ($F \neq\emptyset$) and bounded (optimal objective value is finite) and the recession cone of $f$ and recession cone of $F_R$, the feasible set of continuous relaxation of \textup{\MICP}, have no common non-zero directions of recession, then the continuous relaxation of \textup{\MICP} is bounded.
\end{proposition}
\begin{proof}
Consider $\textup{\MICP}_R$, the continuous relaxation of \textup{\MICP} and $F_R$, the feasible set of $\textup{\MICP}_R$. Let $z_R = \min\set{f({\x})}{\x \in F_R}$ be the optimal objective value of $\textup{\MICP}_R$ (If $\textup{\MICP}_R$ is unbounded then $z_R = -\infty$). Since $z_{R}\leq \zIP$ and $-\infty < \zIP < \infty$, hence the set 
$$
    \bar{F} = \set{\x \in \R^n}{f(\x) \leq \zIP} \neq \emptyset \text{ and } F' = \bar{F} \cap F_R \neq \emptyset.
$$
The recession cone of $F'$ can be represented as,
\begin{align*}
    \textbf{rec}(f)\cap \textbf{rec}(F_R)
\end{align*}
Furthermore, since the recession cone of $f$ and recession cone of $F_R$ have no common non-zero directions of recession, $\textbf{rec}(f)\cap \textbf{rec}(F_R) = \{\bz\}$. This yields that $F'$ is compact. Since $F'$ is compact and non-empty, and $f$ is continuous, $f$ attains a minimum over $F'$. The result follows.
\end{proof}
\begin{lemma}\label{lem:14}
If $f$ is $\mu$-strongly convex and \textup{\MICP} is feasible and bounded then the continuous relaxation of \textup{\MICP} is bounded.
\end{lemma}
\begin{proof}
Let $F \neq \emptyset$ be the feasible set of \MICP. Since $f$ is $\mu$-strongly convex, it follows that for $\x \in \R^n$
\begin{align}
    f(\x) \geq \frac{1}{2}\,\mu\norm{\x-\boldsymbol{\alpha}}^2+\nabla{f(\boldsymbol{\alpha})}^\top(\x-\boldsymbol{\alpha})+ f(\boldsymbol{\alpha}) \text{\hspace{5mm} for some } \boldsymbol{\alpha} \in \R^{n} \label{ineq:strongconv}
\end{align}
Consider the set $F' = \set{\x \in \R^n}{\zIP \geq f(\x)}$. It follows from \cref{ineq:strongconv}, for $\x \in F'$,
\begin{align*}
    \zIP \geq \frac{1}{2}\,\mu\norm{\x-\boldsymbol{\alpha}}^2+ \nabla{f(\boldsymbol{\alpha})}^\top(\x-\boldsymbol{\alpha}) + f(\boldsymbol{\alpha}) \text{\hspace{5mm} for some } \boldsymbol{\alpha} \in \R^{n}
\end{align*}
Consider the set $\displaystyle F'' = \set{\x \in \R^n}{\zIP \geq \frac{1}{2}\,\mu\norm{\x-\boldsymbol{\alpha}}^2+\nabla{f(\boldsymbol{\alpha})}^\top(\x-\boldsymbol{\alpha}) + f(\boldsymbol{\alpha})}$.
It follows that $F' \subseteq F''$. The recession cone of $F''$ is, $$\textbf{rec}(F'') \define \set{\x \in \R^n}{\mathbb{I}\x = \bz, (\nabla{f(\boldsymbol{\alpha}) - \mu\boldsymbol{\alpha})}^\top\x \leq \bz}$$ 
where $\mathbb{I}$ is the $n \times n$ identity matrix \cite{bertsekas2003convex}. Indeed $\textbf{rec}(F'') = \{\bz\}$. This implies that $F''$ is compact and consequently $F'$ is compact. The result follows.
%Hence from \cref{prop:13.1} we get the relaxed set solution exists.
\end{proof}
\begin{lemma}\label{lem:10}
If $\rec(f)$, recession cone of $f$ and the recession cone of $F_R$, the feasible set of continuous relaxation of \textup{\MICP}, have no common non-trivial directions of recession, then for all $0 < \rho < \infty$, the set
$$
S \define \set{\x \in \R^n}{f(\x) + \blambda_{\A}^\top( {\b}- \A\x)+\rho\,\psi({\b}- \A\x) \leq \zIP,\,\E\x \leq \f}
$$
is compact, where $\psi$ is a level-bounding function.
\end{lemma}
\begin{proof}
Consider the Lagrangian function for the continuous relaxation of MICP, i.e. for $\x \in \R^n$
\begin{equation*}
    \mathcal{L}(\x,\blambda_{\A},\blambda_{\E}) = f(\x) - \blambda_\A^{\top}(\A\,\x-\b)-\blambda_\E^{\top}(\f-\E\x)
\end{equation*}
As strong duality holds, the first order necessary (KKT) conditions for $\mathcal{L}(\x,\blambda_{\A},\blambda_{\E})$ can be characterized as
\begin{align}\label{eq:4a}
\begin{split}
    \nabla{f(\x)} = \A^{\top}\blambda_\A - \E^{\top}{\blambda_\E}\\
    \blambda_\E^{\top}(\f - \E\x) = 0\\
    \blambda_\A^{\top}(\A\,\x - \b) = 0\\
    \blambda_\E \geq 0
\end{split}
\end{align}
Observe that $\x_R$, the solution to the continuous relaxation of (\textbf{MICP}) satisfies the system of equations \cref{eq:4a}.\\
Consider $\x \in S$, i.e. $f(\x) + \blambda_{\A}^\top( {\b}- \A\x)+\rho\,\psi({\b}- \A\x) \leq \zIP$. Using convexity of $f$ we obtain,
\begin{align}
    f(\x_R)+ \nabla{f(\x_R)}(\x-\x_R) + \blambda_{\A}^\top( {\b}- \A\x) + \rho\,\psi({\b}- \A\x) \leq \zIP\nonumber\\
    z_R+ \nabla{f(\x_R)}(\x-\x_R) + \blambda_{\A}^\top( {\b}- \A\x) + \rho\,\psi({\b}- \A\x) \leq \zIP \label{ineq:1a}
\end{align}
It follows from strong duality that,
\begin{align*}
\nabla{f(\x_R)}(\x-\x_R) = \blambda_{\A}^{\top} \A(\x - \x_R) -\blambda_\E^{\top} \E(\x-\x_R)\\
\nabla{f(\x_R)}(\x-\x_R) = \blambda_{\A}^{\top} (\A\,\x - \b) -  \blambda_\E^{\top} (\E\x -  \f )
\end{align*}
Substituting in \cref{ineq:1a}, we obtain
\begin{align}
z_R+ \blambda_{\A}^{\top} (\A\,\x - \b) +  \blambda_\E^{\top} ( \f - \E\x ) + \blambda_{\A}^\top( {\b}- \A\x) +\rho\,\psi({\b}- \A\x) & \leq \zIP\nonumber\\
\blambda_\E^{\top} ( \f - \E\x ) +\rho\,\psi({\b}- \A\x) & \leq \zIP - z_R\nonumber\\
\rho\,\psi({\b}- \A\x) & \leq \zIP - z_R \label{ineq:2a}\\
\psi({\b} - \A\x) & \leq \frac{\zIP-z_{R}}{\rho}\nonumber
\end{align}
where \cref{ineq:2a} follows from the fact that $\blambda_\E \geq \bz$ and $\x \in S$. Since $\psi$ is a level-bounded function, there exists a positive $\kappa_\rho$ such that
\begin{align}
     \norm{\b - \A\x}_{\infty} \leq \kappa_{\rho} < \infty\hspace{1cm} \forall\,\,\x \in S.\label{ineq:SS'1}
\end{align}
Additionally, since $f(\x) + \blambda_{\A}^\top({\b} - \A\x) + \rho \psi({\b} - \A\x) \leq \zIP$, $\forall$ $\x \in S$, $\rho > 0$ and $\psi(\u) > 0 \text{ for } \u \neq \bz$
\begin{align}
    f(\x) & \leq \zIP - \blambda_{\A}^\top({\b} - \A\x)\nonumber\\
    & \leq \zIP + \norm{\blambda_{\A}}\norm{\b - \A\x}\nonumber\\
    & \leq \zIP + \sqrt{n}\norm{\blambda_{\A}}\norm{\b - \A\x}_{\infty}\nonumber\\
    & \leq \zIP + \sqrt{n}\norm{\blambda_{\A}}\kappa_{\rho}\label{ineq:SS'2}
\end{align}
Consider the set 
$$
    S' = \set{\x \in \R^n}{f(\x) \leq \zIP + \sqrt{n}\norm{\blambda_{\A}}\kappa_{\rho},\,\norm{\b - \A\x}_{\infty} \leq \kappa_{\rho},\,\E\x \leq \f}.\
$$ 
Inequalities \cref{ineq:SS'1} and \cref{ineq:SS'2} yield that $S \subseteq S'$. Observe that $S'$ can be expressed as,
\begin{align*}
    S' = \set{\x\in \R^n}{\begin{Bmatrix}
f(\x) \leq \zIP + \sqrt{n}\norm{\lambda_{\A}}\kappa_{\rho}\\
{\b} - \kappa_{\rho}\1\leq \A\x \leq {\b} + \kappa_{\rho}\1\\
{\E}\x \leq {\f}
\end{Bmatrix}}
\end{align*}
where ${\1}$ is the vector of ones. We have the recession cone of $S'$ 
\begin{align*}
     \rec(S') = \set{\x\in \R^n}{\begin{Bmatrix}
        \A\x = \bz\\
        {\E}\x \leq \bz
\end{Bmatrix}}\cap \rec(f)
\end{align*}
Alternatively,
\begin{align*}
     \rec(S') = \rec(F_R) \cap \rec(f)
\end{align*}
Since $\rec(F_R) \cap \rec(f) = \{\bz\}$. It follows that $S'$ is a compact set, and consequently $S$ is a compact set.
\end{proof}
\begin{proposition}\label{prop:11}
If $\rec(f)$, recession cone of $f$ and the recession cone of $F_R$, the feasible set of continuous relaxation of \textup{\MICP}, have no common non-trivial directions of recession, then $\lim_{\rho\to\infty} z^{LD+}_{\rho} = \zIP $.
\end{proposition}
\begin{proof}
Recall that the augmented Lagrangian relaxation of \MICP{} is defined as
$$z_{\rho}^{LR+}(\blambda) = \min_{\x \in X} f(\x) + {\blambda}^{\top}(\b - \A\,\x) + \rho\, \psi(\b - \A\,\x)$$
From \cref{lem:10} we have
\begin{align*}
    z^{LR+}_{\rho} & = \min_{\substack{\norm{\x}_\infty \leq M\\\x \in X}}\set{f(\x) + \blambda_{\A}^\top(\b - \A\,\x) + \rho\, \psi(\b - \A\,\x)}{\norm{\b-\A\x}_\infty \leq \kappa_{\rho}}\\
    & \geq \min_{\substack{\norm{\x}_\infty \leq M\\\x \in X}}\set{f(\x) + \blambda_{\A}^\top(\b - \A\,\x)}{\norm{\b-\A\x}_\infty \leq \kappa_{\rho}}%\\
    % & = \min_{\substack{\norm{\x}_\infty \leq M\\\x \in X}}\set{f(\x) + \blambda_{\A}^\top(\b - \A\,\x)}{\norm{\b-\A\x}_\infty \leq \kappa_{\rho}}.
\end{align*}
As a consequence, 
\begin{align}
    \liminf_{\rho\to\infty} z^{LR+}_{\rho} \geq \liminf_{\rho\to\infty} \min_{\substack{\norm{\x}_\infty \leq M\\\x \in X}}\set{f(\x) + \blambda_{\A}^\top(\b - \A\,\x)}{\norm{\b-\A\x}_\infty \leq \kappa_{\rho}}.\label{ineq:liminf}
\end{align}
Define,
\begin{align}
    \theta(\u) & = \min_{\substack{\norm{\x}_\infty \leq M\\\x \in X}}\set{f(\x) + \blambda_{\A}^\top(\b - \A\,\x)}{\A\x = \b + \u} \label{def:theta}
\end{align}
and for all $i \in [k]$,
\begin{align*}
    \Theta(\u, \x^{(i)}_I) & = \min_{\substack{\norm{\x_C}_\infty \leq M\\\left(\x_C,\x^{(i)}_I\right) \in X}}\set{f(\x) + \blambda_{\A}^\top(\b - \A\,\x)}{\A_C\x_C = \b + \u - \A_I\x^{(i)}_I}\\
    \Phi'(\u, \x^{(i)}_I) & = \min_{\substack{\norm{\x_C}_\infty \leq M\\\left(\x_C,\x^{(i)}_I\right) \in X}}\set{f(\x)}{\A_C\x_C = \b + \u - \A_I\x^{(i)}_I},
\end{align*}
and let $\displaystyle\phi'(\u) = \min_{i \in [k]}\Theta(\u, \x^{(i)}_I)$. Observe that, for $\u \in \calU$ and $i \in [k]$, 
\begin{align*}
    \Theta(\u, \x^{(i)}_I) = \Phi'(\u, \x^{(i)}_I) + \blambda_{\A}^\top\u \text{\hspace{5mm}and\hspace{5mm}} \theta(\u) = \min_{i\in [k]} \Theta(\u, \x^{(i)}_I).
\end{align*}
It follows that $\displaystyle\theta(\u) = \min_{i \in [k]}\Phi'(\u, \x^{(i)}_I) + \blambda_{\A}^\top\u = \phi'(\u) + \blambda_{\A}^\top\u$. Furthermore, as a consequence of \cref{ineq:liminf} and \cref{def:theta},
$$
    \liminf_{\rho\to\infty} z^{LR+}_{\rho} \geq \liminf_{\rho\to\infty} \min\set{\theta(\u)}{\norm{\u}_\infty \leq \kappa_\rho}.
$$
Now, as $\rho \to \infty$, $\dfrac{\zIP - z_{R}}{\rho} \to 0$. The level boundedness of $\psi(\cdot)$ implies that as $\dfrac{\zIP - z_{R}}{\rho} \to 0$, $\kappa_\rho \to 0$. Thus, we have
\begin{align*}
    \liminf_{\rho\to\infty} z^{LR+}_{\rho} & \geq \liminf_{\u\to\bz} \phi'(\u) + \blambda_{\A}^\top\u\\
    & = \liminf_{\u\to\bz} \phi'(\u)\\
    & \geq \min\set{f(\x)}{\A\x = \b,\,\norm{\x}_\infty \leq M, \,\x \in X}\\
    & = \zIP.
\end{align*}
The second last inequality follows from lower semi-continuity of $\phi'(\u)$ at $\u = \bz$. The result follows.
\end{proof}
\begin{corollary}\label{corollary:11.1}
If $f$ is strongly convex or if $F_R$ is compact then then $$\lim_{\rho\to\infty} z^{LD+}_{\rho} = \zIP.$$ 
\end{corollary}
%%%%%%%%%%%%%%%%%%%%%%%%%%%%%%%%%%%%%%%%%%%%%%%%%%%%%%%%%%%%%%%%%%%%%%%%%%%%%
\begin{lemma}\label{lem:15}
Let  $\delta > 0$, and let $\u \in \mathcal{N}_{\delta}(\bz)\cap \ocalU$. Define
$$
    S^\u \define  \set{\x \in X \cap \Hu}{ f(\x) \leq \zIP}.
$$
If \textup{\MICP} is feasible ($F \neq\emptyset$) and bounded and the recession cone of $f$ and recession cone of $F_R$, the feasible set of continuous relaxation of \textup{\MICP}, have no common non-zero directions of recession, then  $S^\u$ is bounded.
\end{lemma}
\begin{proof}
Consider $\u \in \mathcal{N}_{\delta}(\bz)\cap \ocalU$. Since $\norm{\u}_\infty \leq \norm{\u} < \delta$ we have $\norm{\b - \A\x}_{\infty} < \delta$, $\forall\,\x \in X \cap \Hu$. Alternatively, $\forall\,\x \in X \cap \Hu$, $-\1\delta < {\b} - \A\x < \delta\1$. Define, for $\u \in \mathcal{N}_{\delta}(\bz)\cap \ocalU$
\begin{align*}
  S^\u \define  \set{\x \in X \cap \Hu}{ f(\x) \leq \zIP}
\end{align*}
Additionally, consider
\begin{align*}
   S = \set{\x\in X_R}{\begin{Bmatrix}
                                    -{\1}\delta \leq {\b} - \A\x \leq {\1}\delta \\
                                    f(\x) \leq z_{IP}
                        \end{Bmatrix}}
\end{align*}
%Here ${\u}$ is taken over all feasible polyhedra i.e. this set $\{\A\x = {\b} + {\u},\,\, {\E}\x\leq {\f}\,\, ,\x\in \mathbb{Z}^{n_1}\times \mathbb{R}^{n_2}\}$ is feasible.
Indeed, $S^\u \subseteq S$. Furthermore, $\rec(S) = \rec(F_R) \cap \rec(f)$, and since $\rec(F_R) \cap \rec(f) = \{\bz\}$, it follows that $S$ is compact and consequently $S^\u$ is bounded. In particular, $\norm{\x}_\infty \leq M < \infty$, $\forall$ $\x \in S^\u$. %since MICP is feasible and optimally bounded S'' is non-empty.\\
%This implies that one can prescribe a uniform bound on variables of the problem and hence $\norm{\x_I\norm{\leq M$.
\end{proof}
The following result is an immediate implication of \cref{lem:15}.
\begin{corollary}\label{cor:16}
If \textup{\MICP} is feasible ($F \neq\emptyset$) and bounded and the recession cone of $f$ and recession cone of $F_R$, the feasible set of continuous relaxation of \textup{\MICP}, have no common non-zero directions of recession, then for ${\u} \in \mathcal{N}_{\delta}(\bz)\cap \ocalU$, $\phi(\u) = \hat{\phi}(\u)$ i.e.
\begin{optprog*}
    {$\phi(\u) = $ minimize}  & \objective{\hspace{1.3cm}f({\x_I},{\x_C})}\\
    subject to & \hspace{0cm} \A_I\,\x_I + \A_C\x_C = \b + \u\\
    & \E_I\x_I + \E_C\x_C \leq \f \hspace{7mm}\\
    & \norm{\x_I}_{\infty}\leq M \hspace{0.5cm}\\
    % & \x_I\in \Z_{+}^{n_1}\hspace{0.4cm}\\ 
     & (\x_I, \x_C) \in \Z^{n_1}\times\R^{n_2}\hspace{0.3cm}
\end{optprog*}
Furthermore, there exists a $\delta > 0$ and $0 < \Gamma < \infty$ such that for every $ \u\in \mathcal{N}_{\delta}(\bz)\cap \calU$, $\phi(\bz)\leq \phi(\u) + \Gamma\norm{\u}$. Additionally, if $f$ is $\mu$-strongly convex and $L$ smooth, then $\Gamma = \mathcal{O}\left(\frac{L\beta\gamma}{\mu}\right)$ where $\displaystyle \beta \define \max_{B \in \mathcal{B}}\norm{B^{-1}}_F$, $\mathcal{B}$ being the set of all possible invertible submatrices of $[\A^T_C \,\, -\A^T_C \hspace{2mm}-{\E}^T_C]$ and $\gamma$ depends on $f(\bz)$, $\norm{\nabla{f(\bz)}}$ and $f(\bar{\x})$ for any $\bar{\x} \in F$.
\end{corollary}
\begin{proof}
\Cref{lem:15}, along with \cref{lem:8,lem:9} readily implies that $\phi(\u) = \hat{\phi}(\u)$ and there exists a $\delta > 0$ and $0 < \Gamma < \infty$ such that for every $ \u\in \mathcal{N}_{\delta}(\bz)\cap \calU$, $\phi(\bz)\leq \phi(\u) + \Gamma\norm{\u}$.

To see the explicit bound on $\Gamma$, recall the stationarity condition from set of equations \cref{eq:32} at $\u = \bz$.
\begin{align}\label{eq:stationarity}
\begin{split}
    \nabla_{C}{f(\x_I^{(i)},\x_C}) & = \A^T_C(\lambda^{(i)+}_{\A_C} - \lambda^{(i)-}_{\A_C}) - {\E}_C^T\lambda_{{\E}_C}\\
    \lambda_{\E_C}^{\top}(\E\x -\f) &= 0\\
    \lambda^{+}_{\A_C},\lambda^{-}_{\A_C} & \geq \bz\\
     \lambda_{{\E}_C} &\geq \bz
\end{split}
\end{align}

For $(\x_C,\x^{(i)}_I)$ consider the complementary slackness conditions, that is if for some $j\in [m]$, $\E_{j}\x < \f_j$ where $\E_{j}$ is a row of the matrix $\E$ then the corresponding Lagrangian multiplier $\lambda_{\E_j} = 0$. Partition the set $[m]$ into two sets $J = \{J_<,J_=\}$, where $J_<=\set{j\in [m]}{\E_j \x <\f_j}$ and $J_= =\set{j\in [m]}{\E_j \x = \f_j}$.\\

Define $\A_{aug} = [\A^T_C \,\, -\A^T_C \hspace{2mm}-{\E}^{\top}_C]$ and $\A_{aug}^J$ as the sub-matrix which has all the columns of $\A_{aug}$ except the columns $2n_2+J_<$.
 % $\lambda^{(i)+}_{\A_C},\lambda^{(i)-}_{\A_C}$. 
 One possible solution to the system \cref{eq:stationarity} is,
 \begin{align}\label{sol:lambda}
    \begin{pmatrix}
     \lambda^{+}_{\A_C}\\
      \lambda^{-}_{\A_C}\\
      \lambda_{{\E}^{J_{=}}_C}
    \end{pmatrix} = B^{-1} \nabla_C{f(\x^{(i)}_I,\x_C)}
\end{align}
where $B$ is a basis matrix of $\A_{aug}^J$.

Furthermore, let $\lambda_{\A_C} = \lambda^{+}_{\A_C} - \lambda^{-}_{\A_C}$, this implies that 
\begin{align}
    \norm{\lambda_{\A_C}} \leq \norm{\lambda^{+}_{\A_C}} + \norm{\lambda^{-}_{\A_C}} \label{ineq:lambda}    
\end{align}
Let $\x^*$ be an optimal solution to \MICP, i.e. $f(\x^*) = \zIP$. Since $f$ is $\mu$-strongly convex, we have

\begin{align*}
    \frac{\mu}{2} \norm{\x^*}^2 + \nabla{f({\bz)}}^\top\x^*  + f(\bz) & \leq z_{IP}\\
    \frac{\mu}{2} \norm{\x^*}^2 + \nabla{f({\bz)}}^\top\x^* - z_{IP} + f(\bz)& \leq 0\\
    \frac{\mu}{2} \norm{\x^*}^2 - \norm{\nabla{f({\bz)}}}\,\norm{\x^*} - \left(z_{IP} - f(\bz)\right) & \leq 0
\end{align*}
It follows that,
$$
    \norm{\x^*}\leq \frac{\norm{\nabla{f(\bz)}} + \sqrt{\norm{\nabla{f(\bz})}^2 +2\,\mu\cdot \left(z_{IP} - f(\bz)\right)}}{\mu}
$$
Since $f$ is L-smooth,
\begin{align*}
    \norm{\nabla{f(\x^*)} - \nabla{f(\bz)}} &\leq L\norm{\x^*}\\
    \norm{\nabla{f(\x^*)}} &\leq L\norm{\x^*} + \norm{\nabla{f(\bz)}}
\end{align*}
From \cref{sol:lambda} and \cref{ineq:lambda} we obtain
\begin{align*}
    \norm{\lambda_{\A_C}} & \leq 2\norm{B^{-1}}_F \norm{\nabla{f(\x^*)}}\\
    & \leq 2\norm{B^{-1}}_F\left(L\norm{\x^*} + \norm{\nabla{f(\bz)}}\right)\\
    & \leq 2\beta\left(L\norm{\x^*} + \norm{\nabla{f(\bz)}}\right)
\end{align*}
where $\beta \define \max_{B \in \mathcal{B}}\norm{B^{-1}}_F$ and $\mathcal{B}$ is the set of all possible invertible submatrices of $[\A^T_C \,\, -\A^T_C \hspace{2mm}-{\E}^T_C]$.

\begin{align}
   \norm{\lambda_{\A_C}} & \leq \dfrac{2\beta\left(L\left(\norm{\nabla{f(\bz)}} + \sqrt{\norm{\nabla{f(\bz})}^2 +2\mu\,\left(z_{IP} - f(\bz)\right)}\right) + \mu\norm{\nabla{f(\bz)}}\right)}{\mu}\nonumber\\
    & \leq  \dfrac{2\beta L}{\mu}\left(2\norm{\nabla{f(\bz)}} + \sqrt{\norm{\nabla{f(\bz)}}^2 + 2\mu\,\left(f(\bar{\x}) - f(\bz)\right)}\right)\label{ineq:thm3}
\end{align}
where the second inequality follows as $\zIP \leq f(\bar{\x}),$ $\forall\,\bar{\x} \in F$ and $L > \mu$. Combining \cref{ineq:thm3} with the construction of $\Gamma$ in \cref{lem:9}, yields the result.
\end{proof}
\cref{thm:3} follows.
\subsection{A Special Case: Pure Integer Convex Programs}
\begin{proposition}\label{prop:7}
Consider the following pure integer convex program \textbf{\textup{(PICP)}},
\begin{optprog*}
    {minimize}  & \objective{\hspace{0.4cm}f(\x)}\\
    \textbf{\textup{(PICP)}}\qquad subject to & \A\,\x = \b \hspace{0.6cm}\\
    & \E\x \leq \f \hspace{7mm}\\
     & \x\in \Z^{n}\hspace{0.4cm}
\end{optprog*}
If the continuous relaxation of \textup{\textbf{(PICP)}} is feasible and bounded then $\exists$ $\rho < \infty$ such that $z^{LD+}_{\rho} = \zIP$.
\end{proposition} 
\begin{proof}
Consider the value function of \textbf{\textup{(PICP)}},
\begin{optprog*}
    {$\hat{\phi}(\u) = $ minimize}  & \objective{\hspace{0.6cm}f(\x)}\\
    $\textbf{\textup{(PICP)}}_{\u}$\qquad subject to & \hspace{0cm} \A\,\x = \b + \u\\
    & \E\x \leq \f \hspace{7mm}\\
    & \x\in \Z^{n}\hspace{0.4cm}
\end{optprog*}
Let $U_{\rho} = \set{\u\in U}{\hat\phi(\u) + \rho\norm{\A\x - \b} \leq z_{IP}}$.
%From \cref{prop:2} and 
From \cref{cor:3} it follows that $\norm{\u}$ is bounded for all $\u \in U_{\rho}$. In particular, we have
\begin{equation*}
    \A\,\x\,= \b + \u\;\; \text{where}\;\norm{\u} \leq \frac{\zIP -  z_{R}}{\rho - \blambda_\A}
\end{equation*}
Since $\A$ and $\b$ are rational, we can assume, without loss of generality, that $\A$ and $\b$ are integral. It follows that there exists an integral solution to the equation $\A\,\x\,= \b +\u$ only if $\u$ is integral. In particular, for $0 < \norm{\u} < 1$ there is no integral solution to the system of equations $\A\,\x\,= \b + \u$. It follows from \cref{prop:2} and \cref{cor:4} that the system of equations $\A\,\x\,= \b +\u$ does not have an integral solution for $\u \in U_\rho$ for any $\rho$ satisfying
\begin{align*}
    \frac{\zIP - z_{R}}{\rho - \norm{\blambda_{\A}}} < 1\\
    (\zIP - z_{R}) < \rho - \norm{\blambda_{\A}} \\
    \rho > (\zIP - z_{R}) + \norm{\blambda_{\A}} = \rho^*
\end{align*}
Thus for $\rho^* < \rho < \infty$ there exists a $\delta_{\rho} > 0$ such that $\forall\, \u \in \mathcal{N}_{\delta_{\rho}}(\bz)\setminus \{\bz\}$, $\textbf{(PICP)}_{\u}$ is infeasible. Furthermore for $\u \in \calU \setminus \mathcal{N}_{\delta_{\rho}}(\bz)$, $\hat{\phi}(\u) + \rho\norm{\u} > \hat{\phi}(\bz)$. The result follows.
\end{proof}

%%%%%%%%%%%%%%%%%%%%%%%%%
%\newpage
\bibliographystyle{siamplain}
\bibliography{Exact_ALD_for_MICP.bib}
\end{document}